\documentclass{article}
\usepackage[T2A]{fontenc}
\usepackage[utf8]{inputenc}
\usepackage[english,russian]{babel}
\usepackage[tbtags]{amsmath}
\usepackage{amsfonts,amssymb,mathrsfs,amscd}

\usepackage{multirow}
\usepackage{subcaption}
\usepackage{float}
\usepackage[nice]{nicefrac}
\usepackage{url}

\usepackage[hyper]{msb-a}

\makeatletter
\gdef\No{{\select@language{russian}\textnumero}}
\makeatother

\JournalName{МАТЕМАТИЧЕСКИЙ СБОРНИК}
\numberwithin{equation}{section}
\theoremstyle{plain}
\newtheorem{theorem}{Теорема}
\newtheorem{lemma}{Лемма}[section]
\newtheorem{propos}{Предложение}
\theoremstyle{definition}
\newtheorem{definition}{Определение}

\newtheorem{remark}{Замечание}

\def\RR{\mathbb R}

\def\R{\mathbb R}

\usepackage{mathtools}
\mathtoolsset{showonlyrefs}
\usepackage{algorithm2e}
\usepackage{algorithmic}
\usepackage{dsfont}

\newcommand{\argmax}{\mathop{\arg\!\max}}
\newcommand{\argmin}{\mathop{\arg\!\min}}

\makeatletter

\def\env@sqcases{
  \let\@ifnextchar\new@ifnextchar
  \left\lbrack
  \def\arraystretch{1.2}
  \array{@{}l@{\quad}l@{}}
}
\makeatother

\newtheorem{re:theorem}{Теорема}
\newtheorem{th:corollary}{Следствие}[theorem]
\newtheorem{re:th:corollary}{Следствие}[re:theorem]
\newtheorem{lm:corollary}{Следствие}[lemma]
\newtheorem{assumption}{Предположение}
\newtheorem{approach}{Подход}

\usepackage[dvipsnames]{xcolor}

\usepackage{makecell}
\usepackage{float}
\usepackage{booktabs,siunitx}

\makeatletter
\gdef\No{{\select@language{russian}\textnumero}}
\makeatother

\begin{document}

\title{Решение сильно выпукло-вогнутых композитных седловых задач с небольшой размерностью одной из групп переменных}
\author[E.\,L.~Gladin]{Е.\,Л.~Гладин}
\address{141701 Московская область, г. Долгопрудный, Институтский переулок, 9, Московский физико-технический институт (национальный исследовательский университет), Россия\\127051 Москва, Большой Каретный пер., 19, Ин-т проблем передачи информации РАН, Россия\\121205 Москва, Большой бульвар, 30с1, Сколковский институт науки и технологий, Россия}
\email{gladin.el@phystech.edu}

\author[I.\,A.~Kuruzov]{И.\,А.~Курузов}
\address{Московский физико-технический институт (национальный исследовательский университет)}
\email{kuruzov.ia@phystech.edu}

\author[D.\,A.~Pasechnyuk]{Д.\,А.~Пасечнюк}
\address{Московский физико-технический институт (национальный исследовательский университет)}
\email{pasechnyuk2004@gmail.com}

\author[F.\,S.~Stonyakin]{Ф.\,С.~Стонякин}
\address{295007 Симферополь, пр-т Акад. Вернадского, 4, Крымский федеральный ун-т, Россия\\Московский физико-технический институт (национальный исследовательский университет)}
\email{fedyor@mail.ru}

\author[M.\,S.~Alkousa]{М.\,С.~Алкуса}
\address{Московский физико-технический институт (национальный исследовательский университет)\\101000 Москва, ул. Мясницкая, 18, НИУ ВШЭ, Россия}
\email{mohammad.alkousa@phystech.edu}

\author[A.\,V.~Gasnikov]{А.\,В.~Гасников}
\address{Московский физико-технический институт (национальный исследовательский университет)\\109004 Москва, ул. А. Солженицына, 25, Центр доверенного искусственного интеллекта Института системного программирования им. В.П. Иванникова РАН, Россия\\
Ин-т проблем передачи информации РАН}
\email{gasnikov.av@mipt.ru}

\date{25.11.2021}
\udk{519.853.62}

\maketitle
\begin{fulltext}

\begin{abstract}
Статья посвящена разработке алгоритмических методов, гарантирующих эффективные оценки сложности для сильно выпукло-вогнутых седловых задач в случае, когда одна из групп переменных имеет большую размерность, а другая~--- достаточно малую (до сотни). Предлагаемая методика основана на сведении задач такого типа к задаче минимизации выпуклого (максимизации вогнутого) функционала по одной из переменных, для которого возможно найти приближённое значение градиента в произвольной точке с необходимой точностью с помощью вспомогательной оптимизационной подзадачи по другой переменной. При этом для маломерных задач предлагается использовать методы эллипсодидов и Вайды, а для многомерных~--- ускоренные градиентные методы с неточной информацией о градиенте или субградиенте. Для случая очень малой размерности задачи одной из групп переменных (до 5) на гиперкубе достаточно эффективным будет иной предлагаемый подход к сильно выпукло-вогнутым седловым задачам на базе нового варианта многомерного аналога метода Ю.\,Е.\,Нестерова на квадрате (многомерная дихотомия) с возможностью использования неточных значений градиента целевого функционала.

Библиография: 35 названий.
\end{abstract}


\markright{}

\footnotetext[0]{
Исследования в пп. 2.2, 2.3 и параграфе 3 выполнены при поддержке Министерства науки и высшего образования Российской Федерации (госзадание) № 075-00337-20-03, номер проекта 0714-2020-0005. Исследования Ф.С. Стонякина в п. 2.1 выполнены при поддержке программы стратегического академического лидерства <<Приоритет - 2030>>, соглашение 075-02-2021-1316 от 30.09.2021.
}

\section{Введение} \label{section_1}
Седловые задачи весьма актуальны, поскольку возникают в реальных проблемах машинного обучения, компьютерной графики, теории игр, а также теории оптимального транспорта. Ввиду важности таких задач известно немало работ, посвящённых различным алгоритмам их решения и теоретическим результатам о скорости их сходимости (сложности) \cite{alkousa2019accelerated,paper:Azizian_2019,paper:Gasnikov_MIPT_2016,paper:Dvinskikh_meta_2020,paper:Hien_Zhao_Haskell,book:Lan2019,paper:Jordan_near_optimal_SPP_2020,usmanova2021fast}.

В данной статье рассматриваются выпукло-вогнутые седловые задачи вида
\begin{equation}\label{problem:min_max00}
    \min_{x \in Q_{x}}\max_{y \in Q_y} \left\{\widehat{S}(x, y) := r(x) + F(x, y) - h(y) \right\},
\end{equation}
где $Q_x\subseteq \mathbb{R}^n, Q_y\subseteq \mathbb{R}^m$~--- непустые выпуклые компактные множества, $r: Q_x \to \mathbb{R}$ и $h: Q_y \to \mathbb{R}$ есть $\mu_{x}$-сильно выпуклая и $\mu_y$-сильно выпуклая функции соответственно. Функционал $F: Q_x \times Q_y \to \mathbb{R}$ выпуклый по $x$ и вогнутый по $y$ и задан в некоторой окрестности множества $Q_x \times Q_y$. Если задача не сильно выпукла (случай $\mu_x=0$ или $\mu_y=0$), то её можно свести к сильно выпуклой применением техники регуляризации (замечание 4.1 в \cite{gasnikov2018book}).

Класс задач \eqref{problem:min_max00} уже некоторое время назад достаточно подробно изучен в билинейном случае, то есть когда $F(x,y)= \langle Ax,y \rangle$ для некоторого линейного оператора $A$ (см., например, обзор \cite{book:Lan2019}). Также известны работы, нацеленные на обобщение известных в билинейном случае результатов на общую ситуацию \cite{paper:Cox_Nemirovsky,paper:Gasnikov_MIPT_2016,paper:Hien_Zhao_Haskell,paper:Nesterov_Excessive}. 

В статье \cite{paper:Zhang_lower_SPP_2019} рассматривалась постановка, когда $Q_x \equiv \RR^n$, $Q_y \equiv \RR^m$, причём для произвольных $x$ и $y$ $\widehat{S}(x, y) = F(x, y)$~--- $\mu_x$-сильно выпуклая по $x$, $\mu_y$-сильно вогнутая по $y$ и $(L_{xx}, L_{xy}, L_{yy})$-гладкая функция. 
Последнее означает, что для любого фиксированного $x$ отображения $\nabla_y F(x, \cdot)$ и $\nabla_x F(x, \cdot)$ являются липшицевыми с некоторыми неотрицательным константами $L_{yy}$ и $ L_{xy}$, а для любого фиксированного $y$ отображения $\nabla_x F(\cdot, y)$ и $\nabla_y F(\cdot, y)$~--- с константами $L_{xx}$ и $L_{xy}$.
В \cite{paper:Zhang_lower_SPP_2019} для выделенного класса задач обоснована нижняя оценка сложности вида $$N(\varepsilon) = \Omega\left(\sqrt{\frac{L_{xx}}{\mu_{x}}+\frac{L_{xy}^{2}}{\mu_x \mu_y}+\frac{L_{yy}}{\mu_y}} \ln\left(\frac{1}{\varepsilon}\right)  \right),$$
где $N(\varepsilon) = \Omega \left( f(\varepsilon) \right) $ означает, что существуют $C>0$ и $\varepsilon_0>0$ такие, что $|N(\varepsilon)| > C|f(\varepsilon)|\ \forall \varepsilon < \varepsilon_0$. В \cite{alkousa2019accelerated} был представлен подход на основе ускоренных методов с оценкой сложности, наиболее близкой к оптимальной на тот момент. Далее, были предприняты попытки получить оптимальный алгоритм \cite{paper:Dvinskikh_meta_2020,paper:Jordan_near_optimal_SPP_2020}. Так, в \cite{paper:Yuanhao_improved} предложен метод с верхней оценкой на количество итераций $\widetilde{O}\left(\sqrt{\frac{L_{xx}}{\mu_x}+\frac{ L \cdot L_{xy}}{\mu_x \mu_y}+\frac{L_{yy}}{\mu_y}} \right)$, где $L=\max \left\{L_{xx}, L_{xy}, L_{yy}\right\}$ (обозначение $\widetilde{O}(\cdot)$ означает $O(\cdot)$ с точностью до логарифмического по $\varepsilon^{-1}$ множителя в степени $1$ или $2$). Таким образом, вопрос о почти оптимальном алгоритме для сильно выпукло-вогнутой седловой задачи большой размерности с гладкой целевой функцией был решён.

В таблице~\ref{review_known_results_for_SPP} приведены наилучшие известные на данный момент результаты
(см. \cite{alkousa2019accelerated,paper:Azizian_2019,paper:Gasnikov_MIPT_2016,paper:Dvinskikh_meta_2020,paper:Hien_Zhao_Haskell,book:Lan2019,paper:Jordan_near_optimal_SPP_2020} и ссылки в них) об оценках сложности решения задачи \eqref{problem:min_max00}.

Для каждого из случаев $\varepsilon$-решение задачи \eqref{problem:min_max00} может быть достигнуто за указанные в первом столбце $\widetilde{O}(\cdot)$ вычислений величины из второго столбца. Константа Липшица $\nabla F$ (градиента по $x$ и по $y$) обозначается как $L_F$.
Далее, когда говорится, что функция $r: Q_x \to \mathbb{R}$ проксимально дружественна, то имеется в виду возможность явно решить задачу вида
\begin{equation}\label{friendly_r}
\min _{x \in Q_{x}}\left\{\langle c_{1}, x\rangle + r(x)+ c_{2}\|x\|_{2}^{2}\right\},\quad c_1 \in Q_x,\ c_2 > 0.
\end{equation}
Аналогичным образом определяется и проксимальная дружественность $h: Q_y \to \mathbb{R}$ для задач вида
\begin{equation}\label{friendly_h}
\min_{y \in Q_y}\left\{\langle c_{3}, y\rangle+ h(y)+ c_{4}\|y\|_{2}^{2}\right\},\quad c_3 \in Q_y,\ c_4 > 0.
\end{equation}
\begin{table}[htp]
\centering
\tabcolsep=0.11cm
\caption{Наилучшие известные 
результаты о сложности методов для задач \eqref{problem:min_max00}. 
}
\label{review_known_results_for_SPP}
\begin{tabular}{|c|c|}
\hline
\multicolumn{2}{|c|}{Случай (1): обе функции $r$ и $h$ проксимально дружественны.}                                                                                   \\ \hline
$\widetilde{O}\left(\frac{L_{F}}{\sqrt{\mu_x \mu_y}} \right)$  & вычислений \eqref{friendly_r}, $\nabla_{x}F(x,y)$ и \eqref{friendly_h}, $\nabla_{y}F(x,y)$ \\ \hline
\multicolumn{2}{|c|}{Случай (2): $r$~--- $L_x$-гладкая и не проксимально дружественная функция.}                                                                    \\ \hline
$\widetilde{O}\left(\sqrt{\frac{L_{x} L_F}{\mu_x \mu_y}}\right)$ & вычислений $\nabla r(x)$                                                                   \\ \hline
$\widetilde{O}\left(\frac{L_{F}}{\sqrt{\mu_x \mu_y}} \right)$ & вычислений $\nabla_{x}F(x,y)$ и \eqref{friendly_h}, $\nabla_{y}F(x,y)$                     \\ \hline
\multicolumn{2}{|c|}{Случай (3): $h$~--- $L_y$-гладкая и не проксимально дружественная функция.}                                                                    \\ \hline
$\widetilde{O}\left(\sqrt{\frac{L_{y} L_F}{\mu_x \mu_y}}\right)$ & вычислений $\nabla h(y)$                                                                   \\ \hline
$\widetilde{O}\left(\frac{L_{F}}{\sqrt{\mu_x \mu_y}} \right)$  & вычислений \eqref{friendly_r}, $\nabla_{x}F(x,y)$  и $\nabla_{y}F(x,y)$                    \\ \hline
\multicolumn{2}{|c|}{\begin{tabular}[c]{@{}c@{}}Случай (4): $r$ и $h$~--- $L_x$- и $L_y$-гладкие\\ не проксимально дружественные функции.\end{tabular}}    \\ \hline
$\widetilde{O}\left(\sqrt{\frac{L_x L_F}{\mu_x \mu_y}}\right)$  & вычислений $\nabla r(x)$                                                                   \\ \hline
$\widetilde{O}\left(\sqrt{\frac{L_y L_F}{\mu_x \mu_y}}\right)$ & вычислений $\nabla h(y)$                                                                   \\ \hline
$\widetilde{O}\left(\frac{L_{F}}{\sqrt{\mu_x \mu_y}} \right)$  & вычислений $\nabla_{x}F(x,y)$ и $\nabla_{y}F(x,y)$                                         \\ \hline
\end{tabular}
\end{table}

Отметим, что седловую задачу можно сводить и к вариационному неравенству (далее~--- ВН) с монотонным оператором. Напомним, что оператор $G: \text{dom}\, G \rightarrow \mathbb{R}^k$, заданный на выпуклом множестве $\text{dom}\, G \subseteq \mathbb{R}^k$, называется монотонным, если
$$\langle G(z) - G(z'),\ z - z' \rangle \geqslant 0\quad \forall z, z' \in Q,$$
где $Q$~--- выпуклое компактное множество с непустой внутренностью и $\text{int}\, Q \subseteq \text{dom}\, G$. Решением вариационного неравенства называется точка $z_* \in \text{dom}\, G \cap Q$, удовлетворяющая соотношению
\begin{equation}\label{VIP}
    \langle G(z_*),\ z - z_* \rangle \geqslant 0\quad \forall z \in Q.
\end{equation}
Выпукло-вогнутые седловые задачи \eqref{problem:min_max0} с дифференцируемой функцией $S(\cdot, \cdot)$ при $r$ и $h$, равных тождественно $0$, сводятся к ВН с оператором $G(x, y) = [\nabla_x S(x, y), -\nabla_y S(x, y)]^\top$, который является монотонным в силу выпуклости $S$ по $x$ и вогнутости по $y$. Такое ВН можно решать, например, с помощью метода эллипсоидов из \cite{nemirovski2010accuracy}, что приводит к скорости сходимости $O \left( \exp \left\{-\frac{N}{2 d(d+1)}\right\} \right)$, где $d = n + m$~--- размерность задачи, $N$~--- число итераций. Такой подход к решению задачи \eqref{problem:min_max0} не требует ни гладкости, ни сильной выпуклости (вогнутости) целевой функции $S(\cdot, \cdot)$ и может считаться вполне эффективным в случае малой размерности задачи $d = n + m$. Более того, базируясь на методах типа метода центров тяжести оценку $O \left( \exp \left\{-\frac{N}{2 d(d+1)}\right\} \right)$ можно улучшить до оценки $O \left( \exp \left\{-\frac{N}{O(d)}\right\} \right)$, см. лекцию 5 \cite{nemirovski1995information}.
Отметим, что первый из указанных выше подходов к задаче \eqref{problem:min_max00} достаточно эффективен в случае большой размерности обеих переменных задачи \eqref{problem:min_max00}, а второй~--- в случае малой размерности переменных задач. 


В настоящей же статье рассматривается ситуация, когда одна из групп переменных ($x$ или $y$) имеет большую размерность, а другая~--- небольшую (несколько десятков). Задача \eqref{problem:min_max00} представляется в виде задачи минимизации по «внешней» переменной выпуклой функции, информация о которой (значение функции, градиента) доступна только с некоторой точностью. Эта точность, в свою очередь, регулируется с помощью вспомогательной оптимизационной подзадачи по «внутренней» переменной. Соответственно, в зависимости от размерности (малая или большая) внешней переменной $x$ естественно выделить два таких подхода. Если размерность внешней переменной $x$ велика (раздел~\ref{1st_approach}), то для решения задачи~\eqref{problem:min_max00} предлагается использовать ускоренный градиентный метод с неточным оракулом, а для вспомогательной максимизационной подзадачи~--- метод секущей гиперплоскости (метод эллипсоидов или метод Вайды). Если же размерность внешней переменной $x$ мала, то для~\eqref{problem:min_max00} используются предложенные в настоящей работе вариации методов секущей гиперплоскости (методы эллипсоидов и Вайды) с использованием на итерациях неточного аналога градиента целевой функции ($\delta$-субградиент или градиент с $\delta$-аддитивной неточностью), а для внутренней оптимизационной подзадачи~--- ускоренные градиентные методы. С целью вывода оценок достаточного количества итераций (обращений к подпрограмме нахождения (суб)градиента функционала или его <<неточного аналога>>) для достижения нужного качества решения седловой задачи ~\eqref{problem:min_max00} при указанном подходе получены важные теоретические результаты, описывающие влияние параметра $\delta$ на качество точки выхода метода эллипсоидов или метода Вайды, которые вместо градиента целевой функции используют на итерациях $\delta$-субградиент или $\delta$-аддитивно неточный градиент. Отметим при этом, что метод Вайды по сравнению с методом эллипсоидов приводит к лучшей оценке достаточного для достижения нужного качества приближённого решения количества итераций, но итерация метода эллипсоидов менее затратна. Далее, детально рассматривается ситуация, когда размерность одной из групп переменных седловой задачи ~\eqref{problem:min_max00} очень мала, а допустимое множество значений этой переменной есть гиперкуб. В этом случае вместо метода эллипсоидов можно использовать аналог метода дихотомии, что может оказаться в случае очень малой размерности (до 5) выгоднее методов секущей гиперплоскости. Точнее говоря, в работе предлагается метод минимизации выпуклой дифференцируемой функции с липшицевым градиентом на многомерном гиперкубе для малых размерностей пространства, который является аналогом метода Ю.\,Е.\,Нестерова минимизации выпуклой липшицевой функции двух переменных на квадрате с фиксированной стороной (см.~\cite{gasnikov2018book,Ston_Pas}). Далее, условимся называть этот метод многомерной дихотомией. Идея метода~--- деление квадрата на меньшие части и постепенное их удаление так, чтобы в оставшейся достаточно малой области все значения целевой функции были достаточно близки к оптимальному. При этом метод заключается в решении вспомогательных задач одномерной минимизации вдоль разделяющих отрезков и не предполагает вычисления точного значения градиента целевого функционала (то есть метод можно считать неполноградиентным). Данный метод в двумерном случае на квадрате рассматривался в работе~\cite{Ston_Pas}. В настоящей же статье предложен новый вариант критерия остановки для вспомогательных подзадач, а также аналог метода Ю.\,Е.\,Нестерова для произвольной размерности. Полученные результаты применимы и к седловым задачам с небольшой размерностью одной из групп переменных на гиперкубе. Приведены оценки сложности такого подхода для сильно выпукло-вогнутых седловых задач вида \eqref{problem:min_max00} с достаточно гладкими функционалами в случае, если маломерная задача решается методом дихотомии с неточным градиентом, и на каждой итерации решается многомерная вспомогательная задача с использованием быстрого градиентного метода для всех вспомогательных задач. Найденные оценки скорости сходимости представляются приемлемыми в случае, если размерность одной из групп переменных достаточно мала (до 5). 

Для сравнения предложенных подходов к задаче \eqref{problem:min_max00} между собой, а также с некоторыми известными аналогами выполнены вычислительные эксперименты для некоторых типов лагранжевых седловых задач к задачам минимизации сильно выпуклых функционалов при наличии небольшого количества выпуклых функциональных ограничений-неравенств. В частности, проведены численные эксперименты для задачи, двойственной к задаче LogSumExp с линейными ограничениями (приложения задач такого типа описаны, например, в работе \cite{paper:Azizian_2019}). Нами проведено сравнение скорости работы в случае, когда маломерная задача решается быстрым градиентным методом с $(\delta,L)$-оракулом и методом эллипсоидов с $\delta$-субградиентом. Проведённые вычислительные эксперименты показали, что использование рассмотренных методов секущей гиперплоскости для задач небольшой размерности приводит к более удачной работе с достаточно высокой требуемой точностью по сравнению с методами, использующими только градиентные подходы. Кроме этого, в ряде случаев метод многомерной дихотомии показал себя эффективнее методов секущей гиперплоскости, что указывает на целесообразность и такого предлагаемого нами подхода. Для случая, когда маломерная подзадача решается методом эллипсоидов, был дополнительно проведён эксперимент по сравнению методик учёта неточностей при использовании аддитивно зашумлённого градиента (раздел~\ref{sec:exp-additive-noise}): предложенный подход со значением неточности, изменяющимся вместе с диаметром текущего эллипсоида, позволяет применяемому методу быстрее достигать условия останова и, следовательно, заданной точности. Проведён также эксперимент по сравнению эффективности предложенных подходов с известными альтернативами применительно к задаче проектирования точки на множество, заданное набором (небольшого числа) гладких ограничений (раздел~\ref{sec:exp-proj}). Сравнение показало, что подход с применением метода эллипсоидов для маломерной задачи при использовании предложенного условия останова и новых оценок на достаточное число итераций для внутреннего метода оказывается существенно эффективнее аналогичного подхода, предложенного в работе \cite{usmanova2021fast}. Отметим, что для применяемых в данной работе к подзадачам небольшой размерности методов секущей гиперплоскости или многомерного аналога дихотомии сильная выпуклость важна лишь для теоретических оценок. Для реализации таких методов предполагать сильную выпуклость целевой функции не обязательно, и поэтому мы не проводим регуляризацию лагранжевых седловых задач по двойственным переменным.

Работа состоит из введения, заключения и двух основных разделов. В разделе~\ref{LabelSect2} приводятся основные результаты, подходы к рассматриваемой задаче для различных случаев малой размерности внешних и внутренних переменных. В п.~\ref{LabSect2.1} описана общая схема рассуждений, используемых для анализа рассматриваемых в настоящей работе задач, которая связана с рассмотрением семейства вспомогательных подзадач оптимизации. Следующий пункт~\ref{1st_approach} является ключевым и содержит вывод оценок для задачи \eqref{problem:min_max00} в случае относительно небольшой размерности одной из групп переменных на базе использования новых вариаций методов эллипсоидов и Вайды к соответствующим вспомогательным подзадачам. Далее, в п.~\ref{3st_approach} рассмотрен специальный случай очень малой (до 5) размерности одной из групп переменных седловой задачи и описаны оценки сложности подхода к \eqref{problem:min_max00}, основанного на предлагаемом многомерном аналоге дихотомии с аддитивно неточным градиентом. В разделе~\ref{labsect3} приводятся результаты некоторых вычислительных экспериментов и сравнение скорости работы предложенных подходов. Отметим, что полные доказательства некоторых результатов (теорем~\ref{th_ellips},~\ref{th:dich_x},~\ref{InexGradConst},~\ref{CurGrad},~\ref{small} и леммы~\ref{subgradient}) приводятся в приложении к работе.


\section{Основные результаты}\label{LabelSect2}

\subsection{Схема вывода оценок сложности для рассматриваемого класса седловых задач}\label{LabSect2.1}

\subsubsection{Постановка задачи}
\label{general_alg}
Перепишем задачу \eqref{problem:min_max00} в виде
\begin{equation}\label{problem:min_max0}
    \min_{x \in Q_{x}} \left\{ g(x):= r(x) + \max_{y \in Q_y} S(x, y) \right\},
\end{equation}
где
$Q_x \subseteq \mathbb{R}^n, Q_y \subseteq \mathbb{R}^m$~--- выпуклые замкнутые множества, $Q_x$~--- ограничено, $S(x, y) = F(x, y) - h(y)$~--- непрерывная функция из \eqref{problem:min_max00}, сильно выпуклая по $x$ и сильно вогнутая по $y$.
\begin{definition}\label{definepsaccuracy}
Будем называть пару точек $(\tilde{x},\tilde{y}) \in Q_x \times Q_y$ $\varepsilon$-решением задачи \eqref{problem:min_max00} (или \eqref{problem:min_max0}), если
\begin{equation}\label{eqsaddlesolution}
\max\{\|\tilde{x} - x_*\|_2, \, \|\tilde{y} - y_*\|_2\} \leqslant \varepsilon,
\end{equation} 
где $(x_*, y_*)$~--- точное решение сильно выпукло-вогнутой седловой задачи \eqref{problem:min_max00}.
\end{definition}


\begin{remark}
Отметим, что ввиду сильной выпуклости функции $g$ для выполнения неравенства \eqref{eqsaddlesolution} достаточно потребовать находить для задачи \eqref{problem:min_max0} такое $\tilde{x}$, что $g(\tilde{x}) - \min_{x \in Q_{x}} g(x) \leqslant C\varepsilon^2$ (при подходящем выборе константы $C>0$, зависящей от параметра сильной выпуклости $\mu_x$), а также с аналогичной точностью решать вспомогательные подзадачи. Поскольку используем методы с гарантией линейной скорости сходимости, то в итоговые оценки сложности для задач \eqref{problem:min_max00} величины типа $C\varepsilon^2$ будут входить под логарифмами. Это означает, что для вывода асимптотических оценок сложности седловых задач вида \eqref{problem:min_max00} достаточно потребовать нахождения такого $\tilde{x}$, чтобы $g(\tilde{x}) - \min_{x \in Q_{x}} g(x) \leqslant \varepsilon$, что мы и будем использовать в рассуждениях.
\end{remark} 

Будем рассматривать задачу \eqref{problem:min_max0} как композицию внутренней задачи максимизации
\begin{equation}\label{problem:max_S0}
   \widehat{g}(x):= \max_{y \in Q_y} S(x, y)
\end{equation}
и внешней задачи минимизации
\begin{equation}\label{problem:min_g0}
    \min_{x \in Q_x} g(x).
\end{equation}
Итерационный метод для внешней задачи \eqref{problem:min_g0} будет на каждом шаге использовать градиент целевой 
функции, который может быть вычислен с некоторой точностью на основе приближённого решения внутренней задачи \eqref{problem:max_S0}. В связи с этим возникает необходимость в чётких оценках качества выдаваемого методом решения в случае использования на его итерациях неточной информации о градиенте или субградиенте целевой функции. Оказывается, что для седловых задач в качестве подходящего неточного аналога субградиента целевой функции можно рассматривать $\delta$-субградиент (см. определение~\ref{DeltaSubgrad}), $(\delta, L)$-субградиент (см. \eqref{deltaLsubgrad} ниже) \cite{alkousa2019accelerated} или $\delta$-неточный субградиент (см. определение~\ref{DefInexactGrad} ниже).

\begin{definition}\label{DeltaSubgrad}
    Пусть $\delta \geqslant 0$. Вектор $\nu(\hat{x}) \in \mathbb{R}^n$ называется $\delta$-субградиентом выпуклой функции $g : Q_x \to \mathbb{R}$ в точке $\hat{x}$, если $g(x) \geqslant g(\hat{x}) + \langle \nu(\hat{x}), x - \hat{x} \rangle - \delta$ для всякого $x \in Q_x$. Множество $\delta$-субградиентов $g$ в $\hat{x}$ обозначается $\partial_{\delta} g(\hat{x})$.
\end{definition}
\noindent Заметим, что $\delta$-субградиент совпадает с обычным субградиентом при $\delta=0$.

\begin{definition}\label{DefInexactGrad}
    Пусть $\delta \geqslant 0$. Будем называть вектор $\nu(\hat{x}) \in \mathbb{R}^n$ $\delta$-неточным субградиентом выпуклой функции $g : Q_x \to \mathbb{R}$ в точке $\hat{x}$, если для некоторого субградиента $\nabla g(\hat{x}) \in \partial g(\hat{x})$ выполнено $\left\| \nabla g(\hat{x}) - \nu(\hat{x}) \right\|_2 \leq \delta$. Если известно, что функция $g$ дифференцируема в точке $\hat{x}$, то будем говорить, что $\nu(\hat{x})$~--- $\delta$-неточный градиент.
\end{definition}

\subsubsection{Вычисление неточного аналога градиента целевой функции основной подзадачи}

В качестве приближённого субградиента целевой функции $g$ задачи \eqref{problem:min_g0} в точке $x \in Q_x$ 
предлагается использовать субградиент $\nabla r(x) + \nabla_x S(x, \tilde{y})$,
где $\tilde{y}$~--- $\tilde{\varepsilon}$-решение вспомогательной подзадачи \eqref{problem:max_S0} при данном $x$, $\nabla r(x) \in \partial r(x)$ и $\nabla_x S(x, \tilde{y}) \in \partial_x S(x, \tilde{y})$ ~--- произвольные конечные субградиенты в точке $x$ функций $r$ и $S(\cdot, \tilde{y})$ соответственно. Оказывается, что такой неточный субградиент может быть $\delta$-субградиентом целевой функции $g$ в точке $x$, если выбрать точность $\tilde{\varepsilon}$ для вспомогательной подзадачи согласно следующей лемме.
\begin{lemma}\label{lem:ellips}
    (см. также \cite{polyak1983intro}, с. 123--124) Пусть в условиях задачи \eqref{problem:min_max0} для фиксированного $x$ точка $\widetilde{y} \in Q_y$ такова, что $\widehat{g}(x) - S(x, \widetilde{y}) \leqslant \delta$. Тогда $\partial_x S(x, \widetilde{y}) \subseteq \partial_{\delta} (\widehat{g}(x))$.
\end{lemma}

\noindent Приведённая лемма говорит, что для отыскания $\delta$-субградиента функции $g$ достаточно решить задачу максимизации \eqref{problem:max_S0} с точностью $\tilde{\varepsilon} = \delta$.

  
Оказывается \cite{alkousa2019accelerated}, что для седловых задач типа \eqref{problem:min_max0} при соответствующих предположениях и точности решения вспомогательной подзадачи можно гарантировать при подходящем $L>0$ и сколь угодно малом $\delta > 0$ доступность $(\delta, L)$-субградиента $\nabla_{\delta, L} g(x)$ функции $g$ в произвольной точке $x \in Q_x$:
\begin{equation}\label{deltaLsubgrad}
g(x) + \langle \nabla_{\delta, L} g(x), y - x \rangle - \delta \leqslant g(y) \leqslant g(x) +  \langle \nabla_{\delta, L} g(x), y - x \rangle + \frac{L}{2}\|y - x\|_2^2 + \delta.   
\end{equation}
Ясно, что $(\delta, L)$-субградиент $\nabla_{\delta, L} g(x)$ есть некоторый $\delta$-субградиент функции $g$ в точке $x$ с дополнительным условием вида \eqref{deltaLsubgrad}.

Следующий известный результат непосредственно вытекает из аналогичного результата \cite{th2_cite} для известного понятия ($\delta, L$)-оракула и поясняет связь между двумя используемыми далее аналогами градиента ($\delta$-субградиент и $\delta$-неточный субградиент) выпуклой функции $g$, допускающей ($\delta, L$)-липшицев субградиент в каждой точке $x \in Q_x$.
\begin{theorem}\label{lem:boundary}
    Пусть $g: Q_x \to \mathbb{R}$~--- выпуклая функция, $\nu(x) = \nabla_{\delta, L} g(x)$~--- её $(\delta, L)$-субградиент в точке $x \in \operatorname{int} Q_x$. В таком случае, если $\rho (x, \partial Q_x)$~--- евклидово расстояние от точки $x$ до границы множества $Q_x$, при $\rho (x, \partial Q_x) \geqslant 2\sqrt{\frac{\delta}{L}}$ для всякого субградиента $\nabla g(x)$ верно неравенство
$$\| \nu(x) - \nabla g(x) \|_2 \leqslant 2\sqrt{\delta L}.$$
\end{theorem}
\noindent На базе леммы~\ref{lem:ellips} и теоремы~\ref{lem:boundary} можно заключить, что при достаточно малом $\delta > 0$ для отыскания $\delta$-неточного субградиента функции $g$ в точке $x \in \operatorname{int} Q_x$ достаточно решить задачу максимизации \eqref{problem:max_S0} с точностью $\tilde{\varepsilon}= \frac{\delta^2}{4L}$.

Ещё один способ нахождения $\delta$-неточного субградиента (в данном случае уже $\delta$-неточного градиента, поскольку дополнительно предполагается дифференцируемость) может быть использован при следующем дополнительном предположении.
\begin{assumption}\label{asssumption:delta_vs_L}
Функция $S(\cdot, y)$ дифференцируема для всех $y \in Q_y$ и удовлетворяет для некоторого $L_{xy} \geq 0$ условию
    \begin{equation}\label{Lxy0}
        \|\nabla_x S (x, y) - \nabla_x S(x, y')\|_2\leqslant L_{xy} \|y - y' \|_2\quad \forall x \in Q_x,\, y,y' \in Q_y.
    \end{equation}
\end{assumption}
\begin{lemma}\label{lem:delta_vs_L}
    Пусть в условиях задачи \eqref{problem:min_max0} и предположения~\ref{asssumption:delta_vs_L} для любого фиксированного $x\in Q_x$ точка $\widetilde{y} \in Q_y$ такова, что $\widehat{g}(x) - S(x, \widetilde{y}) \leqslant \tilde{\varepsilon}$. Тогда $\widehat{g}$  дифференцируема в точке $x$ и выполняется неравенство
    \begin{equation*}
        \|\nabla_x S (x, \widetilde{y}) - \nabla \widehat{g}(x)\|_2 \leqslant L_{xy} \sqrt{\frac{2 \tilde{\varepsilon}}{\mu_y}}.
    \end{equation*}
\end{lemma}
Таким образом, для отыскания $\delta$-неточного градиента функции $g$ достаточно решить задачу максимизации \eqref{problem:max_S0} с точностью $\tilde{\varepsilon}= \frac{\mu_y}{2 L_{xy}^2} \delta^2$, если $L_{xy} >0$, и с любой конечной точностью $\tilde{\varepsilon}$, если $L_{xy} =0$.

Приведём ещё одну лемму о связи двух аналогов градиента ($\delta$-субградиента и $\delta$-неточного субградиента), которые далее встречаются в статье.
\begin{lemma}\label{delta_12}
Пусть $g: Q \to \mathbb{R}$~--- выпуклая функция на выпуклом множестве $Q$, тогда
\begin{enumerate}
    \item если множество $Q$ ограничено, то $\delta_1$-неточный субградиент функции $g$ является её $\delta_2$-субградиентом с $$\delta_2 = \delta_1 \operatorname{diam} Q, \text{ где } \operatorname{diam} Q = \sup_{x, x' \in Q} \|x-x'\|_2;$$
    \item если функция $g$ $\mu$-сильно выпуклая, то её $\delta_1$-неточный субградиент является $\delta_2$-субградиентом с $\delta_2 = \frac{\delta_1^2}{2\mu}$.
\end{enumerate}
\end{lemma}

Методику исследований и результаты статьи можно условно подразделить на две части, первая из которых (пункт~\ref{1st_approach}) связана c использованием методов секущей гиперплоскости (методы эллипсоидов и Вайды) для подзадач небольшой размерности, а вторая (пункт~\ref{3st_approach}) использует авторский многомерный вариант дихотомии с адаптивными правилами остановки для подзадач малой размерности (примерно до 5). Для результатов, основанных на использовании для основных подзадач методов секущей гиперплоскости, приведённые выше утверждения позволяют обосновать возможность использовать на итерациях как $\delta$-субградиент, так и $\delta$-неточный субградиент целевой функции $g$. Оценки сложности для седловых задач при этом будут асимптотически совпадать. Что же касается второй части результатов, связанных с использованием многомерной дихотомии, то для них существенно используется допущение о гладкости функции и $\delta$-неточного градиента (именно уже градиента, поскольку в этой части рассматриваются только гладкие задачи).

\subsubsection{Общая схема (алгоритм) подхода к выделенному классу задач}
В настоящем пункте приводится алгоритм для минимаксной задачи \eqref{problem:min_max0}, а в 
последующих разделах рассматриваются конкретные примеры методов, используемых в общем алгоритме, и соответствующие оценки сложности.

\begin{algorithm}
	\caption{Алгоритм для минимаксной задачи \eqref{problem:min_max0}.}
	\label{alg:general}
	\begin{algorithmic}[1]
		\REQUIRE Метод $\mathcal{M}_1$ для решения задачи \eqref{problem:min_g0} с использованием $\delta$-субградиента или $\delta$-неточного градиента, число шагов этого метода $N > 0$, метод $\mathcal{M}_2$ для решения задачи \eqref{problem:max_S0}, точность $\tilde{\varepsilon}$ её решения, начальное приближение $(x^0, y^0)$
		\FOR{$k=0,\, \dots, \, N-1$}
		    \STATE Решить задачу \eqref{problem:max_S0} при фиксированном $x = x^k$ с точностью $\tilde{\varepsilon}$ методом $\mathcal{M}_2$, стартуя из $y^k$:
		    $$y^{k+1} := \mathcal{M}_2(x^k, y^k, \tilde{\varepsilon})$$
		    \STATE Положить $\nu^{k+1} := \nabla r(x^k) + \nabla_x S(x^k, y^{k+1}) \in \partial r(x^k) + \partial_x S(x^k, y^{k+1})$
		    \STATE Сделать один шаг метода $\mathcal{M}_1$ из точки $x^k$, используя приближённый градиент $\nu^{k+1}$:
		    $$x^{k+1} := \text{step } (\mathcal{M}_1, x^k, \nu^{k+1})$$
		\ENDFOR
		\ENSURE $x^N$.
	\end{algorithmic}
\end{algorithm}
Сложность алгоритма~\ref{alg:general} будем определять согласно следующему очевидному принципу.
\begin{propos}\label{th:general_alg}
    Пусть метод $\mathcal{M}_1$ для решения задачи \eqref{problem:min_g0} с использованием $\delta$-субградиента или $\delta$-неточного градиента находит $\varepsilon$-решение не более чем за $N_1(\varepsilon, \delta)$ шагов\footnote{Для того, чтобы конкретный метод $\mathcal{M}_1$ мог гарантированно обеспечить точность $\varepsilon$ за конечное число шагов, может потребоваться, чтобы $\delta$ было достаточно малым по сравнению с $\varepsilon$ (например, $\delta < \varepsilon$). Предложение предполагает, что такое условие выполнено.}, и пусть метод $\mathcal{M}_2$ для решения задачи \eqref{problem:max_S0} находит $\tilde{\varepsilon}$-решение не более чем за $N_2(\tilde{\varepsilon})$ шагов. Если точность $\delta$ оракула для задачи \eqref{problem:min_g0} зависит от $\tilde{\varepsilon}$ как $\delta(\tilde{\varepsilon})$, то алгоритм~\ref{alg:general} находит $\varepsilon$-решение задачи \eqref{problem:min_max0} после $N_1(\varepsilon, \delta(\tilde{\varepsilon}))$ вычислений $\nabla_x S$ и $N_1(\varepsilon, \delta(\tilde{\varepsilon})) \cdot N_2(\tilde{\varepsilon})$ вычислений $\nabla_y S$.
\end{propos}

\subsection{Методы секущей гиперплоскости с использованием неточных аналогов субградиента и их приложения к оценкам сложности для седловых задач с небольшой размерностью одной из групп переменных}
\label{1st_approach}

Приведём конкретные методы, которые могут быть использованы в алгоритме~\ref{alg:general} в случае, когда размерность внешней или внутренней переменной относительно мала (не более ста), а целевая функция имеет композитную структуру. 

Начнём с постановки задачи и краткого описания методов, а также полученных оценок сложности в случае малой размерности основной подзадачи (иными словами, внешней переменной).

Пусть для задачи \eqref{problem:min_max0} выполнено следующее
\begin{assumption}\label{assum_for_vaidya}
Множество $Q_x$ имеет непустую внутренность, размерность $n$ относительно мала (не более ста), $Q_y \equiv \mathbb{R}^m$, функция $S$ имеет вид
\begin{equation}\label{composite_function}
    S(x, y) := F(x,y) - h(y),
\end{equation}
где $\mu_y$-сильно выпуклая функция $h$ непрерывна, выпукло-вогнутая функция $F$ дифференцируема по $y$ и удовлетворяет для некоторого $L_{yy} \geq 0$ условию
\begin{equation}\label{smooth_F_1_2}
    \|\nabla_y F(x, y)-\nabla_y F(x, y')\|_2 \leqslant L_{yy}\|y-y'\|_2\quad \forall x \in Q_x,\, y,y' \in Q_y.
\end{equation}
Пусть также выполнено одно из условий:
\begin{enumerate}
    \renewcommand{\theenumi}{\alph{enumi}}
    \item \label{assum_prox} $h$ является проксимально-дружественной, то есть в явном виде решается задача
    \begin{equation}\label{prox_h}
        \min_{y \in Q_{y}}\left\{\langle c_{1}, y\rangle + h(y)+ c_{2}\|y\|_{2}^{2}\right\},\quad c_1 \in Q_y, c_2 > 0;
    \end{equation}
    \item \label{assum_smooth_h} $h$ имеет $L_h$-липшицев градиент.
\end{enumerate}
\end{assumption}

\begin{theorem}\label{complexity_theorem}
$\varepsilon$-решение задачи \eqref{problem:min_max0} в предположении~\ref{assum_for_vaidya} может быть достигнуто за
\begin{equation}
    O \left( n \ln \frac{n}{\varepsilon} \right) \text{ вычислений } \nabla_x F,\, \nabla r \text{ и}
\end{equation}
\begin{itemize}
    \item в предположении~\ref{assum_for_vaidya}.\ref{assum_prox}~---
    \begin{equation}
        O \left( n \sqrt{\frac{L_{yy}}{\mu_y}} \ln \frac{n}{\varepsilon} \ln \frac{1}{\varepsilon} \right) \text{ вычислений } \nabla_y F \text{ и решений } \eqref{prox_h};
    \end{equation}
    \item в предположении~\ref{assum_for_vaidya}.\ref{assum_smooth_h}~---
    \begin{gather}
        O \left( n \sqrt{\frac{L_{yy}}{\mu_y}} \ln \frac{n}{\varepsilon} \ln \frac{1}{\varepsilon} \right) \text{ вычислений } \nabla_y F,\\
        O \left( n \sqrt{\frac{L_h}{\mu_y}} \ln \frac{n}{\varepsilon} \ln \frac{1}{\varepsilon} \right) \text{ вычислений } \nabla h.
    \end{gather}
\end{itemize}
\end{theorem}

Далее, приводятся методы, применимые к вспомогательным подзадачам из алгоритма~\ref{alg:general}, а также теоретические результаты об оценках их скорости сходимости.

\subsubsection{Методы секущей плоскости
с использованием $\delta$-субградиентов}
Рассмотрим задачу вида
\begin{equation}\label{problem:ell}
    \min_{x \in Q} g(x),
\end{equation}
где $Q \subset \R^n$~--- выпуклое компактное множество, которое содержится в некотором евклидовом шаре радиуса $\mathcal{R}$ и включает некоторый евклидов шар радиуса $\rho>0$, $g$~--- непрерывная выпуклая функция, число $B>0$ таково, что $|g(x) - g(x')| \leqslant B\ \forall x, x' \in Q$.

Мы предлагаем обобщение метода эллипсоидов (алгоритм~\ref{alg:ellipsoid}) для задачи \eqref{problem:ell}, на итерациях которого используется $\delta$-субградиент целевой функции.
\begin{algorithm}
	\caption{Метод эллипсоидов с $\delta$-субградиентом для задачи \eqref{problem:ell}.}
	\label{alg:ellipsoid}
	\begin{algorithmic}[1]
		\REQUIRE Число итераций $N > 0$, $\delta \geqslant 0$, шар $\mathcal{B}_{\mathcal{R}} 	\supseteq Q$, его центр $c$ и радиус $\mathcal{R}$.
		\STATE $\mathcal{E}_0 := \mathcal{B}_{\mathcal{R}},\quad H_0 := \mathcal{R}^2 I_n,\quad c_0 := c$.
		\FOR{$k=0,\, \dots, \, N-1$}
		    \IF {$c_k \in Q$}
		        \STATE $w_k := w \in \partial_{\delta} g(c_k)$, 
		        \IF {$w_k = 0$}
		            \RETURN $c_k$, 
		        \ENDIF
		    \ELSE
		        \STATE $w_k := w$, где $w \neq 0$ таков, что $Q \subset \{ x \in \mathcal{E}_k: \langle w, x-c_k \rangle \leqslant 0 \}.$
		    \ENDIF
		    \STATE $c_{k+1} := c_k - \frac{1}{n+1}\frac{H_k w_k}{\sqrt{w_k^T H_k w_k}}$, \\
		    $H_{k+1} := \frac{n^2}{n^2-1} \left( H_k - \frac{2}{n+1}\frac{H_k w_k w_k^T H_k}{w_k^T H_k w_k} \right)$, \\
		    $\mathcal{E}_{k+1} := \{x: (x-c_{k+1})^T H_{k+1}^{-1} (x-c_{k+1}) \leqslant 1 \}$,
		\ENDFOR
		\ENSURE $x^N = \arg\min\limits_{x \in \{c_0, \ldots, c_N \} \cap Q } g(x)$.
	\end{algorithmic}
\end{algorithm}
\begin{theorem}[Оценка качества приближённого решения для метода эллипсоидов с использованием $\delta$-субградиентов]\label{th_ellips}
    После $N \geqslant 2n^2 \ln \left(\frac{\mathcal{R}}{\rho}\right)$ итераций алгоритм~\ref{alg:ellipsoid} для задачи \eqref{problem:ell} возвращает такую точку $x^N \in Q$, что
    \begin{equation}\label{th_ellips_1}
        g(x^N) - \min_{x \in Q} g(x) \leqslant \frac{B \mathcal{R}}{\rho} \exp \left(-\frac{N}{2n^2} \right)+\delta.
    \end{equation}
\end{theorem}

\begin{th:corollary}

Если дополнительно к условиям теоремы~\ref{th_ellips} допустить, что $g$ $\mu_x$-сильно выпукла, то точка выхода $x^N$ алгоритма~\ref{alg:ellipsoid} удовлетворяет неравенству
    \begin{equation}\label{th_ellips_2}
        \| x^N - x_* \|_2^2 \leqslant \frac{2}{\mu_x} \left( \frac{B \mathcal{R}}{\rho} \exp \left(-\frac{N}{2n^2} \right)+\delta \right),
    \end{equation}
    где $x_*$~--- искомая точка минимума $g$.
\end{th:corollary}

\begin{remark}
Условие $\mu_x$-сильной выпуклости $g$ и оценка \eqref{th_ellips_2} существенны для обоснования достижимости требуемого качества решения седловой задачи \eqref{problem:min_max00} согласно определению~\ref{definepsaccuracy}.
\end{remark}


\begin{remark}\label{inexactnesses_ellipsoid}
Отметим, что метод эллипсоидов с $\delta$-субградиентом может быть использован и в том случае, когда у нас есть доступ к $\delta$-неточному градиенту вместо точного или $\delta$-субградиента, см. лемму~\ref{delta_12}.
\end{remark}

Теперь напомним метод секущей плоскости, который был предложен Вайдой \cite{vaidya1989new} для решения задачи \eqref{problem:ell}.
Сначала введём необходимые обозначения. Для матрицы $A$ и вектора $b$ будем рассматривать вспомогательный ограниченный $n$-мерный многогранник $P(A, b)$ вида
\begin{equation}\label{polytope}
    P(A, b) = \{x \in \R^n: \, Ax\geq b\}, \text{ где } A \in \R^{m\times n},\, b \in \R^m,
\end{equation}
где неравенство $Ax\geq b$ понимается как покомпонентное (каждая координата вектора $Ax$ не меньше соответствующей координаты вектора $b$).

Для множества $P(A, b)$ можно ввести следующий логарифмический барьер 
$$
L(x; A,b) = -\sum_{i=1}^{m} \ln \left(a_{i}^{\top} x-b_{i}\right),\quad x \in \operatorname{int} P(A,b),
$$
где $a_{i}^{\top}$~--- $i$-я строка матрицы $A$, $\, \operatorname{int} P(A,b)$~--- внутренность $P(A,b)$. Гессиан $H$  функции $L$ равен
\begin{equation}\label{hess}
    H(x; A,b) =\sum_{i=1}^{m} \frac{a_{i} a_{i}^{\top}}{\left(a_{i}^{\top} x-b_{i}\right)^{2}},\quad x \in \operatorname{int} P(A,b).
\end{equation}
Матрица $H(x; A,b)$ положительно определена для всех $x \in \operatorname{int} P(A,b)$. Также для множества $P(A, b)$ можно ввести волюметрический барьер (volumetric barrier) вида
\begin{equation}\label{vol_center}
    V(x; A,b) = \frac{1}{2} \ln \left(\operatorname{det}H(x; A,b)\right),\quad x \in \operatorname{int} P(A,b),
\end{equation}
где $\operatorname{det}H(x; A,b)$ обозначает определитель $H(x; A,b)$. Обозначим за $\sigma_{i}(x; A,b)$ величины
\begin{equation}\label{sigmas}
    \sigma_{i}(x; A,b)=\frac{a_{i}^{\top} \left(H(x; A,b)\right)^{-1} a_{i}}{\left(a_{i}^{\top} x-b_{i}\right)^{2}},\quad x \in \operatorname{int} P(A,b), \quad 1 \leq i \leq m.
\end{equation}
Волюметрическим центром множества $P(A, b)$ называют точку минимума волюметрического барьера
\begin{equation}\label{vol_center}
    x_c = \argmin_{x \in \operatorname{int} P(A, b)} V(x; A,b). 
\end{equation}

Волюметрический барьер $V$ является самосогласованной функцией, поэтому может быть эффективно минимизирован методом Ньютона. Подробный теоретический анализ для метода Вайды можно найти в статье \cite{vaidya1989new} и книге \cite{bubeck2015convex}. Ранее в \cite{gladin2021solving} доказано, что в методе Вайды можно использовать $\delta$-субградиент вместо точного субградиента. Ниже приводится вариант метода с использованием $\delta$-субградиента (алгоритм~\ref{alg:vaidya}). Этот алгоритм образует последовательность пар $\left(A_k, b_k\right) \in \R^{m_k\times n}\times \R^{m_k}$ таких, что соответствующие многогранники содержат искомое решение задачи. В качестве исходного многогранника, задаваемого парой $\left(A_0, b_0\right)$, можно выбрать, например, симплекс
\begin{equation}
    P_0=\Bigl\{x \in \mathbb{R}^n: x_{j} \geqslant-\mathcal{R}, j=\overline{1,n},\ \sum_{j=1}^{n} x_{j} \leqslant n \mathcal{R} \Bigr\} \supseteq \mathcal{B}_{\mathcal{R}} \supseteq \mathcal{X},
\end{equation}
то есть
\begin{equation}\label{A_and_b}
    b_0 = - \mathcal{R} \left[\begin{array}{c}
        \mathbf{1}_n \\
        n
    \end{array}\right],\quad A_0 = \left[\begin{array}{c}
        I_n \\ 
        -\mathbf{1}_n^{\top}
    \end{array}\right],
\end{equation}
где $I_n$ обозначает единичную матрицу размера $n \times n$, $\mathbf{1}_n$ обозначает вектор из единиц $(1, \dots, 1)^{\top} \in \mathbb{R}^n$. В таком случае $m_0$ будет равно $n+1$.


\begin{algorithm}[h!]
	\caption{Метод Вайды с использованием $\delta$-субградиентов для задач \eqref{problem:ell}.}
	\label{alg:vaidya}
	\begin{algorithmic}[1]
		\REQUIRE Число итераций $N > 0$, $\delta \geqslant 0$, пара $(A_0, b_0)$
		(см. \eqref{A_and_b}), $m_0:=n+1$, параметры алгоритма $\eta \leqslant 10^{-4}$, $\gamma \leqslant 10^{-3} \cdot \eta$.
		\FOR{$k=0,\, \dots, \, N-1$}
		    \STATE Найти приближённый волюметрический центр, см. 
		    \eqref{vol_center}. 
		    \STATE Вычислить $H_k^{-1} := \left( H(x_k; A_k,b_k) \right)^{-1}$ и $\displaystyle \left\{ \sigma_{i}(x_k; A_k,b_k) \right\}_{i=1}^{m_k}$ по формулам \eqref{hess} и \eqref{sigmas},
		    \STATE $\displaystyle i_k := \arg \min_{1 \leqslant i \leqslant m_k} \sigma_{i}(x_k; A_k,b_k)$
		    \IF {$\sigma_{i_k}(x_k; A_k,b_k) < \gamma$}
		        \STATE Получить $\left(A_{k+1}, b_{k+1}\right)$ исключением $i_k$-й строки из $\left(A_k, b_k\right)$, 
		        \STATE $m_{k+1} := m_k - 1.$
		    \ELSE
		        \STATE $c_k \in -\partial_\delta g(x_k)$,
		        \STATE Найти $\beta_k \in \mathbb{R}$, удовлетворяющий $c_k^\top x_k \geq \beta_k$, из уравнения
		        $$\frac{c_k^\top H_k^{-1} c_k}{(c_k^\top x_k - \beta_k)^2} = \frac{1}{2} \sqrt{\eta \gamma},$$
		        \STATE $A_{k+1} := \begin{pmatrix}A_k\\c_k^{\top}\end{pmatrix},\;\;b_{k+1} := \begin{pmatrix}b_k\\\beta_k\end{pmatrix},\;\;m_{k+1} = m_k + 1$.
		    \ENDIF
		\ENDFOR
		\ENSURE $x_N = \arg\min\limits_{x \in \{x_0, ..., x_{N-1}\}} g(x)$.
	\end{algorithmic}
\end{algorithm}


\begin{theorem}{\cite{gladin2021solving}}\label{th:vaidya}
    После $N \geq \frac{2n}{\gamma} \ln \left( \frac{n^{1.5} \mathcal{R}}{\gamma \rho} \right) + \frac{1}{\gamma} \ln \pi$ итераций метод Вайды с $\delta$-субградиентом для задачи \eqref{problem:ell} возвращает такую точку $x^N$, что
\begin{equation}\label{accuracy_vaidya}
    g(x^N) - \min_{x \in Q} g(x) \leqslant \frac{n^{1.5} B \mathcal{R}}{\gamma \rho} \exp \left( \frac{\ln \pi -\gamma N}{2n} \right) + \delta,
\end{equation}
где $\gamma>0$~--- параметр алгоритма~\ref{alg:vaidya}.
\end{theorem}
\begin{th:corollary}
Если дополнительно к условиям теоремы~\ref{th:vaidya} добавить $\mu_x$-силь\-ную выпуклость $g$, то точка выхода $x^N$ алгоритма~\ref{alg:vaidya} удовлетворяет неравенству
    \begin{equation}\label{th_ellips_21}
        \| x^N - x_* \|_2^2 \leqslant \frac{2}{\mu_x} \left( \frac{n^{1.5} B \mathcal{R}}{\gamma \rho} \exp \left( \frac{\ln \pi -\gamma N}{2n} \right) + \delta \right),
    \end{equation}
где $x_*$~--- точка минимума $g$. 
\end{th:corollary}

\begin{remark}[учёт неточной информации о значении целевой функции]\label{inexact_function}
Отметим, что как метод эллипсоидов, так и метод Вайды используют значения целевой функции при определении выходов алгоритмов ($x^N$). Однако по смыслу рассматриваемой в настоящей статье постановки седловой задачи естественна ситуация, когда значение целевой функции вспомогательной подзадач доступно лишь с некоторой точностью $\tilde{\delta}$. В таком случае выписанные оценки качества выдаваемых методами приближённых решений \eqref{th_ellips_1} и \eqref{accuracy_vaidya} необходимо уточнить, добавив в правые части неравенств дополнительные слагаемые $\tilde{\delta}$. Действительно, если функция $g_{\tilde{\delta}}$ отличается от $g$ не больше чем на $\tilde{\delta}$, то для
$$
\tilde{x}^N := \arg\min\limits_{x \in \{x_0, ..., x_{N-1}\}} g_{\tilde{\delta}}(x)\, \text{ и }\, x^N := \arg\min\limits_{x \in \{x_0, ..., x_{N-1}\}} g(x)
$$
верно неравенство $g(\tilde{x}^N) \leq g(x^N) + \tilde{\delta}$.
\end{remark}

\begin{remark}
    Метод Вайды с использованием $\delta$-субградиентов может быть использован и в том случае, когда методу доступна информация о $\delta$-неточном субградиенте вместо точного или $\delta$-субградиента, см. лемму~\ref{delta_12}.
\end{remark}
\begin{remark}[сравнение результатов о сложности для метода эллипсоидов и метода Вайды]\label{mat_inversion}
    По числу итераций, требуемых для достижения заданной точности решения минимизационной задачи по функции, метод эллипсоидов проигрывает методу Вайды. Действительно, для метода эллипсоидов оценка количества итераций квадратично зависит от размерности пространства, а для метода Вайды эта оценка пропорциональна $n \ln n$. С другой стороны, сложность одной итерации у метода эллипсоидов меньше, чем у метода Вайды. Действительно, на итерации метода Вайды необходимо находить обратную матрицу к квадратной матрице порядка $n$, а в первом~--- достаточно лишь ограничиться операцией умножения матрицы такого размера на вектор.
\end{remark}

\subsubsection{Ускоренные методы для задач композитной оптимизации в пространствах больших размерностей}
Если предыдущий пункт содержит изложение методов и теоретических результатов, применяемых для возникающих при решении основных задач \eqref{problem:min_max00} или \eqref{problem:min_max0} подзадач небольшой размерности, то сейчас будут рассмотрены используемые подходы к вспомогательным подзадачам большой размерности. Точнее говоря, опишем методы для задач выпуклой композитной минимизации вида
\begin{equation}\label{problem:um}
    \min_{y \in \mathbb{R}^m} \left\{ U(y) := u(y) + v(y) \right\},
\end{equation}
где $u$~--- $\mu$-сильно выпуклая функция с $L_u$-липшицевым градиентом, $v$~--- выпуклая функция.

\begin{algorithm} [h!]
\caption{Ускоренный метаалгоритм (УМ) \cite{paper:Dvinskikh_meta_2020} для задачи \eqref{problem:um}.}
\label{alg:um}
	\begin{algorithmic}[1]
		\REQUIRE Число итераций $K \geqslant 1$, начальная точка $z_0$, параметр $H>0$.
		\STATE $A_0 = 0,\quad y_0 = z_0$.
		\FOR{ $k = 0, \ldots, K-1$}
		\STATE $ \lambda_{k+1}  = \displaystyle \frac{1}{2H}, $ 
		\STATE $ a_{k+1} = \displaystyle \frac{\lambda_{k+1}+\sqrt{\lambda_{k+1}^2+4\lambda_{k+1}A_k}}{2},\ A_{k+1} = A_k+a_{k+1}, $
		\STATE $ \widetilde{z}_k = \displaystyle \frac{A_k}{A_{k + 1}}y_k + \frac{a_{k+1}}{A_{k+1}} z_k $,
		\STATE \begin{equation}\label{problem:um_aux}
    		y_{k+1} = \argmin_{y\in \RR^d} \left\{ u(\widetilde{z}_k) +\langle \nabla u(\widetilde{z}_k) ,y - \widetilde{z}_k\rangle + v(y)  +\frac{H}{2}\|y-\widetilde{z}_k\|_2^{2} \right\},
		\end{equation}
		\STATE $z_{k+1} := z_k-a_{k+1} \nabla u(y_{k+1})
 		- a_{k+1}\nabla v(y_{k+1})$,
		\ENDFOR
		\ENSURE УМ($z_0$, $K$) $:=y_{K}$.
	\end{algorithmic}
\end{algorithm}

\begin{algorithm} [h!]
\caption{Рестартованный УМ \cite{paper:Dvinskikh_meta_2020}.}
\label{alg:restart_um}
	\begin{algorithmic}[1]
		\REQUIRE Количество рестартов $K \geqslant 1$, начальная точка $z_0$, параметры $H, \mu>0$.
		\FOR{ $k = 0, \ldots, K-1$}
		\STATE $\displaystyle N_k=\left\lceil \sqrt{ \frac{32H}{\mu}} \right\rceil$,
		\STATE  $z_{k+1} :=$ УМ($z_k$,$N_k$) (алгоритм~\ref{alg:um}),
		\ENDFOR
		\ENSURE $z_{K}$.
	\end{algorithmic}	
\end{algorithm}
\begin{theorem}[Оценка сложности рестартованного УМ \cite{paper:Dvinskikh_meta_2020}]\label{theorem:um}
    Пусть $z_N$~--- выход алгоритма~\ref{alg:restart_um} после $N$ итераций. Тогда если $H \geqslant 2 L_u$, то общее число вычислений \eqref{problem:um_aux} для достижения $\displaystyle U(z_N) - U(y_*) \leqslant \varepsilon$ будет равно
    \begin{equation}\label{compl_rest_um}
        N = O \left( \sqrt{\frac{H}{\mu}} \ln \left(\frac{\mu R_y^2}{\varepsilon}\right) \right),
    \end{equation}
    где $R_y = \| y^0 - y_* \|_2$, а $y_*$~--- точное решение задачи~\eqref{problem:um}.
\end{theorem}

\begin{remark}[разделение оракульных сложностей]\label{slide_um}
Если $v$ имеет $L_v$-лип\-ши\-цев градиент, то можно смотреть на вспомогательную задачу \eqref{problem:um_aux} как на гладкую сильно выпуклую задачу. Для её решения можно также использовать рестартованный УМ ($u^{new} := v,\,v^{new}(y) := u(\widetilde{z}_k) + \langle \nabla u(\widetilde{z}_k), y - \widetilde{z}_k\rangle +\frac{H}{2}\|y-\widetilde{z}_k\|_2^{2},\,\mu^{new} := \frac{H}{2},\,H^{new} := 2 L_v$), подобно тому, как это сделано в \cite{paper:Dvinskikh_meta_2020}.
При условии $L_u \leqslant L_v$ это позволяет получить $\varepsilon$-решение \eqref{problem:um} за
\begin{equation}\label{complexity:um}
    O \left( \sqrt{\frac{H}{\mu}} \ln \left(\frac{\mu R_y^2}{\varepsilon}\right) \right) \text{ вычислений } \nabla u \text{ и } O \left( \sqrt{\frac{L_v}{\mu}} \ln \left(\frac{\mu R_y^2}{\varepsilon}\right) \right) \text{ вычислений } \nabla v.
\end{equation}
\end{remark}

\subsubsection{Оценки сложности алгоритмов для седловых задач с использованием методов эллипсоидов и Вайды для подзадач небольшой размерности}\label{second_x}

Для того, чтобы найти $\varepsilon$-решение задачи \eqref{problem:min_max0} в предположении~\ref{assum_for_vaidya}, предлагается использовать следующий
\begin{approach}[Случай малой размерности $x$]\label{appr:first}
Алгоритм~\ref{alg:general} применяется к задаче \eqref{problem:min_max0} с параметрами $\tilde{\varepsilon}:= \frac{\varepsilon}{2}$, $\mathcal{M}_1$~--- метод Вайды (алгоритм~\ref{alg:vaidya}), $\mathcal{M}_2$~--- рестартованный УМ (алгоритм~\ref{alg:restart_um}).
\end{approach}

Воспользуемся предложением~\ref{th:general_alg}, чтобы установить сложность подхода~\ref{appr:first}.
Для этого требуется,
используя обозначения из формулировки упомянутого предложения,
выписать зависимости $N_1(\varepsilon, \delta)$, $\delta(\tilde{\varepsilon})$ и $N_2(\tilde{\varepsilon})$.
Согласно лемме~\ref{lem:ellips} точность  $\tilde{\varepsilon}= \frac{\varepsilon}{2}$ решения внутренней задачи приводит к точности $\delta$-субградиента, равной $\delta = \frac{\varepsilon}{2}$. Как видно из теоремы~\ref{th:vaidya}, количество итераций метода Вайды составляет
\begin{equation}
    N_1 \left( \varepsilon, \frac{\varepsilon}{2} \right) = \left\lceil \frac{2 n}{\gamma} \ln \left(\frac{2 n^{1.5} B \mathcal{R}}{\gamma \rho \varepsilon} \right) + \frac{\ln \pi}{\gamma} \right\rceil.
\end{equation}
Рестартованный УМ (алгоритм~\ref{alg:restart_um}) применяется к функциям $u(\cdot) := -F(x, \cdot)$ и $v(\cdot) := h(\cdot)$ (если в предположении~\ref{assum_for_vaidya}.\ref{assum_smooth_h} $L_{yy} > L_h$, то нужно поменять $u$ и $v$). Прежде чем выписать его сложность, заметим, что отображение $y^*(x):= \argmax_{y \in Q_y} S(x, y)$ непрерывно в силу непрерывности $S$ и её сильной выпуклости по $y$. Следовательно, множество $\{ y^*(x)\; |\; x \in Q_x \}$ ограничено как образ компактного множества. Обозначим его диаметр через $R_y$.

Если $h$ является проксимально-дружественной \textbf{(случай 1)}, то согласно формуле \eqref{compl_rest_um} $\frac{\varepsilon}{2}$-решение внутренней задачи может быть найдено за
\begin{equation*}
    N_2 \left( \frac{\varepsilon}{2} \right) = O \left( \sqrt{\frac{L_{yy}}{\mu_y}} \ln \left(\frac{\mu_y R_y^2}{\varepsilon}\right) \right) \text{ вычислений } \nabla_y F \text{ и } \eqref{prox_h}.
\end{equation*}
Если же $h$ имеет $L_h$-липшицев градиент \textbf{(случай 2)}, то согласно формулам \eqref{complexity:um} $\frac{\varepsilon}{2}$-решение внутренней задачи может быть найдено за
\begin{equation*}
    N_2^F \left( \frac{\varepsilon}{2} \right) = O \left( \sqrt{\frac{L_{yy}}{\mu_y}} \ln \left(\frac{\mu_y R_y^2}{\varepsilon}\right) \right) \text{ вычислений } \nabla_y F \text{ и }
\end{equation*}
\begin{equation*}
    N_2^h \left( \frac{\varepsilon}{2} \right) = O \left( \sqrt{\frac{L_h}{\mu_y}} \ln \left(\frac{\mu_y R_y^2}{\varepsilon}\right) \right) \text{ вычислений } \nabla h.
\end{equation*}
Из выписанных оценок сложности и предложения~\ref{th:general_alg} следует результат, сформулированный в теореме~\ref{complexity_theorem}. Вместо метода Вайды можно использовать в качестве $\mathcal{M}_1$ метод эллипсоидов. Согласно теореме~\ref{th_ellips} его сложность составляет
\begin{equation}
    N_1 \left( \varepsilon, \frac{\varepsilon}{2} \right) = \left\lceil 2 n^2 \ln \left( \frac{2 B \mathcal{R}}{\rho \varepsilon} \right) \right\rceil.
\end{equation}
В этом случае получаются похожие оценки сложности, а именно, $\varepsilon$-решение задачи \eqref{problem:min_max0} в предположении~\ref{assum_for_vaidya} может быть достигнуто за $O \left( n^2 \ln \frac{1}{\varepsilon} \right)$ вычислений $\nabla_x F,\, \nabla r$ и
\begin{itemize}
    \item в предположении~\ref{assum_for_vaidya}.\ref{assum_prox}~---
    $O \left( n^2 \sqrt{\frac{L_{yy}}{\mu_y}} \ln^2 \frac{1}{\varepsilon} \right)$ вычислений $\nabla_y F$ и решений \eqref{prox_h};
    \item в предположении~\ref{assum_for_vaidya}.\ref{assum_smooth_h}~---
    $O \left( n^2 \sqrt{\frac{L_{yy}}{\mu_y}} \ln^2 \frac{1}{\varepsilon} \right)$ вычислений $\nabla_y F$, и 
    
    $O \left( n^2 \sqrt{\frac{L_h}{\mu_y}} \ln^2 \frac{1}{\varepsilon} \right)$ вычислений $\nabla h$.
\end{itemize}

Отметим, что когда функция $h$ является проксимально-дружественной (предположение~\ref{assum_for_vaidya}.\ref{assum_prox}), то в подходе~\ref{appr:first} можно использовать в качестве метода $\mathcal{M}_2$, например, метод подобных треугольников \cite{gasnikov2018mpt}, который позволяет убрать из предположения~\ref{assum_for_vaidya} требование $Q_y \equiv \mathbb{R}^m$, сохранив те же оценки сложности.

Другой интересный случай возникает, когда $Q_y = [a_1, b_1] \times \dots \times [a_m, b_m]$~--- гиперпрямоугольник и $S$ сепарабельна по $y$, то есть для $y = (y_1, y_2, \ldots, y_m) \in Q_y$ верно $S(x, y) = \sum_{i=1}^m S_i(x, y_i)$, где для любого $x \in Q_x$ функции $S_i(x, y_i)$ аргумента $y_i$ являются непрерывными и унимодальными. Тогда из предположения~\ref{assum_for_vaidya} можно убрать требование о гладкости и сильной вогнутости функции $S$ по $y$ и сводить вспомогательную задачу \eqref{problem:max_S0} к $m$ задачам одномерной максимизации
$$
\max_{y_i \in [a_i, b_i]} S_i(x, y_i),\ i = \overline{1, m}.
$$
Эти задачи можно решать с помощью метода дихотомии (деления отрезка пополам) с точностью $\frac{\varepsilon}{2m}$, что гарантирует точность $\frac{\varepsilon}{2}$ решения вспомогательной задачи \eqref{problem:max_S0}. Такой подход позволяет достигнуть $\varepsilon$-решение задачи \eqref{problem:min_max0} за $O \left( n \ln \frac{n}{\varepsilon} \right)$ вычислений $\nabla_x F$ и $O \left( mn \ln \frac{n}{\varepsilon} \ln \frac{m}{\varepsilon} \right)$ вычислений $S(x, y)$. 

Наконец, рассмотрим случай, когда малую размерность имеют не внешние переменные ($x$), а внутренние ($y$). Пусть
$F$ является выпуклой по $x$ и $\mu_y$-сильно вогнутой по $y$, а также для любых $x \in Q_x$ и $y,\,y' \in Q_y$, верны неравенства:
\begin{equation}\label{smooth_F_1}
    \|\nabla_x F(x, y)-\nabla_x F(x', y)\|_2 \leqslant L_{xx}\|x-x'\|_2,
\end{equation}
\begin{equation}\label{smooth_F_2}
    \|\nabla_x F(x, y)-\nabla_x F(x, y')\|_2 \leqslant L_{xy}\|y-y'\|_2,
\end{equation}
\begin{equation}\label{smooth_F_3}
    \|\nabla_y F(x, y)-\nabla_y F(x', y)\|_2 \leqslant L_{xy}\|x-x'\|_2,
\end{equation}
где  $L_{xx}, L_{xy} < \infty$. Пусть также функция $r$ является $\mu_x$-сильно выпуклой и проксимально-дружественной. 
В такой постановке внешняя задача \eqref{problem:min_g0} (минимизация функции $g$) является многомерной.
Как показано в \cite{alkousa2019accelerated}, возможно использовать для минимизации функции $g$ быстрый градиентный метод с аналогом неточного оракула~--- $(\delta, L)$-моделью целевой функции в произвольной запрашиваемой точке $Q_x$. Поэтому для решения задачи \eqref{problem:min_max0} предлагается использовать следующий
\begin{approach}[малая размерность $y$]\label{appr:second}
Для внешней задачи \eqref{problem:min_g0} применяется быстрый градиентный метод для $(\delta, L)$-оракула для задач сильно выпуклой композитной оптимизации (\cite{alkousa2019accelerated}, алгоритм 4). Для внутренней задачи \eqref{problem:max_S0} используется метод Вайды (алгоритм~\ref{alg:vaidya} с $\delta=0$).
\end{approach}
Можно доказать, что такой подход позволяет достигнуть $\varepsilon$-решение задачи \eqref{problem:min_max0} за
$O \left( \sqrt{\frac{L}{\mu_x}} \ln \frac{1}{\varepsilon} \right)$ вычислений $\nabla_x F$ и решений \eqref{friendly_r}, $O \left( m \sqrt{\frac{L}{\mu_x}} \ln \frac{1}{\varepsilon} \ln \frac{m}{\varepsilon} \right)$ вычислений $\nabla_y F,\, \nabla h$, где $L = L_{xx} + \frac{2L_{xy}^2}{\mu_y}$.

\subsection{Метод многомерной дихотомии для задач оптимизации малой размерности на гиперкубе и её приложения к седловым задачам}\label{3st_approach}

В настоящем подразделе мы рассмотрим выпукло-вогнутую седловую задачу (полагаем одно из композитных слагаемых $r$ в задаче \eqref{problem:min_max00} тождественно равным 0)
\begin{equation}\label{max_min}
    \max_{y \in Q_{y}} \left\{ \min_{x \in Q_{x}} S(x, y) :=  F(x,y) - h(y) \right\}.
\end{equation}

Пусть для задачи \eqref{max_min} выполнено следующее
\begin{assumption}\label{assum_for_dich}
$Q_x\subseteq \mathbb{R}^n$~--- гиперкуб с конечной стороной, $Q_y\subseteq \mathbb{R}^m$~--- непустое выпуклое компактное множество, размерность $n$ очень мала (менее 5), функция $\widehat{S}$ имеет вид
\begin{equation}\label{composite_function}
   \widehat{S}(x, y) = S(x, y) = F(x,y) - h(y),
\end{equation}
где $\mu_y$-сильно выпуклая функция $h$ непрерывна, функционал $F$, заданный в некоторой окрестности множества $Q_x \times Q_y$, является выпуклым по $x$ и вогнутым по $y$. Допустим, что $F$ является достаточно гладким, а именно, для произвольных $x, x' \in Q_x$ и $y, y' \in Q_y$ верны неравенства:
    \begin{equation}\label{smooth_F_11}
        \|\nabla_x F(x, y)-\nabla_x F(x', y)\|_2 \leqslant L_{xx}\|x-x'\|_2,
    \end{equation}
    \begin{equation}\label{smooth_F_22}
        \|\nabla_x F(x, y)-\nabla_x F(x, y')\|_2 \leqslant L_{xy}\|y-y'\|_2,
    \end{equation}
    \begin{equation}\label{smooth_F_33}
        \|\nabla_y F(x, y)-\nabla_y F(x', y)\|_2 \leqslant L_{xy}\|x-x'\|_2,
    \end{equation}
    \begin{equation}\label{smooth_F_44}
        \|\nabla_y F(x, y)-\nabla_y F(x, y')\|_2 \leqslant L_{yy}\|y-y'\|_2,
    \end{equation}
где  $L_{xx}, L_{xy}, L_{yy} <+\infty$.
Пусть, как и ранее (см. предположение~\ref{assum_for_vaidya}), выполнено одно из условий~\ref{assum_prox} (\textbf{случай 1}) или~\ref{assum_smooth_h} (\textbf{случай 2}).
\end{assumption}

Введем для задачи \eqref{max_min} функцию $f$:
\begin{equation}
    \label{f_max_min}
    f:=\max_{y\in Q_y} S(x,y).
\end{equation}
Обозначим диаметр множества $Q_x$ через $R=\max_{x_1,x_2\in Q_x}\|x_1-x_2\|$. Если размер каждой стороны гиперкуба равен $a$, то $R=a\sqrt{n}.$

Как и в предыдущих подпунктах, будем рассматривать подходы к \eqref{max_min} на базе системы вспомогательных минимизационных задач. Однако для решения маломерных (внешних) подзадач на гиперкубе уже будем использовать некоторый аналог метода дихотомии. Поэтому сначала опишем этот подход к решению задач выпуклой минимизации на многомерном гиперкубе с использованием неточного градиента на итерациях.

Рассмотрим задачу оптимизации следующего вида:
\begin{equation}\label{func_f_1}
\min_{x \in Q_x} f(x),
\end{equation}
где функция $f$ является липшицевой с константой $M_f$, имеет липшицев градиент с константой $L_f$ и является $\mu_f$-сильно выпуклой, $Q_x$~--- конечный гиперкуб. В этом разделе статьи будет описан и проанализирован алгоритм~\ref{alg:Dichotomy} (многомерная дихотомия), являющийся аналогом предложенного Ю.\,Е.\,Нестеровым метода двумерной минимизации на квадрате  \cite{gasnikov2018mpt}, упр. 4.2~--- метод многомерной дихотомии на гиперкубе размерности $n\geqslant 2$ (см. алгоритм~\ref{alg:Dichotomy}).

Для алгоритма~\ref{alg:Dichotomy} получен следующий результат.

\begin{theorem}\label{th:dich_x}
Пусть для задачи \eqref{func_f_1} функция $f$ является $M_f$-липшицевой $\mu_f$-сильно выпуклой и имеет $L_f$-липшицев градиент. Для достижения точности $\varepsilon$ по функции точки выхода алгоритма~\ref{alg:Dichotomy} достаточно
\begin{equation}\label{dich_x}
O\left(2^{n^2} \log_2^n\left(\frac{C R}{\varepsilon}\right)\right),\text{ где } C = \max\left(M_f,\frac{4(M_f+2L_fR)}{L_f}, \frac{128L_f^2}{\mu_f}\right)
\end{equation}
обращений к подпрограмме для вычисления $\nu(\textbf{x})$, где $\nu(\textbf{x})$ есть приближение градиента $\nabla f$ такое, что $\|\nabla f(\textbf{x}) - \nu(\textbf{x})\|_2\leqslant \widetilde{\delta}(\textbf{x})$ для всякой текущей точки $\textbf{x}$. Точность $\widetilde{\delta}(x)$ определяется из условия \eqref{Adapt_inexact}, точность решения вспомогательных задач определяется согласно \eqref{ConstEstFunc}.
\end{theorem}

Используя этот результат и методы вспомогательной минимизации, описанные в разделе~\ref{1st_approach}, можно прийти к следующим выводам. Для достижения точности $\varepsilon$ в смысле определения~\ref{definepsaccuracy} для седовой задачи \eqref{max_min} при использовании многомерной дихотомии для решения по малоразмерноей переменной, достаточно:
\begin{equation}\label{dich_x_fgm_y}
\begin{aligned}
O\left(2^{n^2} \log_2^n\left(\frac{C R}{\varepsilon}\right)\right)
\text{ обращений к подпрограмме вычисления $\nabla_x S(x, y)$ и}
\end{aligned}
\end{equation}
\begin{enumerate}
    \item в случае 1 ---
    
        $O \left( 2^{n^2} \sqrt{\frac{L_{yy}}{\mu_y}} \log_2^{n+1} \left(\frac{CR}{\varepsilon}\right) \right)$  вычислений  $\nabla_y F(x, y)$ и решений подзадачи \eqref{prox_h};
    \item в случае 2 ---
    \begin{equation}\label{N_h}
    \begin{aligned}
        & O \left(  \sqrt{\frac{L_h}{\mu_y}} \log_2^{n+1}\left(\frac{CR}{\varepsilon}\right) \right) \text{вычислений } \nabla h(y) \\
        & \text{и } O \left( 2^{n^2} \left( \sqrt{\frac{L_h}{\mu_y}}+\sqrt{\frac{L_{yy}}{\mu_y}} \right) \log_2^{n+1}\left(\frac{CR}{\varepsilon}\right) \right) \text{ вычислений } \nabla_y F(x, y).
    \end{aligned}
    \end{equation}
\end{enumerate}

\subsubsection{Описание метода многомерной дихотомии}\label{theor_dich}

Метод многомерной дихотомии предполагает проведение через центр гиперкуба разделяющей гиперплоскости параллельно одной из его граней и рекурсивное решение вспомогательной задачи оптимизации с некоторой точностью $\widetilde{\varepsilon}$, выбор которой мы обсудим ниже. В точке приближённого решения вычисляется неточный градиент $\nu(x)$, такой, что $\|\nu(x) - \nabla f(x)\|_2 \leqslant \widetilde{\delta}$ для подходящего $\widetilde{\delta}$, и выбирается та его компонента в текущей точке, которая соответствует зафиксированному измерению выбранной гиперплоскости, и в зависимости от её знака выбирается та часть гиперкуба, в которой не лежит неточный градиент. На одной итерации метода описанная процедура выполняется для каждой грани гиперкуба. Итерации метода выполняются для основного гиперкуба до тех пор, пока размер $R$ оставшейся области не станет меньше $\frac{\varepsilon}{M_f}$, что гарантирует сходимость по функции с точностью $\varepsilon$. Условие останова для вспомогательной подзадачи, гарантирующее достижение приемлемой точности исходной задачи, будет подробнее описано ниже (см. теорему~\ref{FullCond}).

\begin{algorithm}[h!]
    \caption{Многомерная дихотомия.}
    \label{alg:Dichotomy}
    \begin{algorithmic}[1]
		\REQUIRE множество $Q = \bigotimes\limits_{k=1}^n[a_k, b_k]$, требуемая точность $\varepsilon$ по функции, процедура вычисления неточного градиента $\nu(\textbf{x})$, возвращающая элемент множества $\left\{\nu(\textbf{x}) \Big | \|\nu(\textbf{x}) - \nabla f(\textbf{x})\|_2\leqslant \widetilde{\delta}\right\}$, начальное приближение $\textbf{x}$, требуемое число итераций $N^*.$
    \IF{$Q=\{x\}$}
        \RETURN $x$
    \ENDIF
    \WHILE{$N \leqslant N^*$}
        \FOR{$i=1, \ldots, n$}
            \STATE $c:=\frac{a_i+b_i}{2}$
            \STATE Фиксация одной из размерностей:
            $$Q_{\text{new}} := \left\{x\in Q\Big| x_i = c\right\}$$
            \STATE Подпрограмма для вычисления неточного градиента во вспомогательной задаче:$$\nu_{\text{new}}(\textbf{x}) := \left(\nu_1(\textbf{x})\;\;\dotsi\;\; \nu_{i-1}(\textbf{x})\;\;\nu_{i+1}(\textbf{x})  \;\;\dotsi \;\;\nu_n(\textbf{x}) \right)$$
            \STATE Рекурсивный вызов многомерной дихотомии для гиперкуба с размерностью, уменьшенной на один, $Q_{\text{new}}$ и с новой требуемой точностью решения $\widetilde{\varepsilon} = \frac{\mu_f}{128 L^2 R^2}\varepsilon^2$ (см. \eqref{ConstEstFunc})
            $$\textbf{x} := \text{Dichotomy}(Q_{\text{new}}, \widetilde{\varepsilon}, \nu_{\text{new}})$$
            \STATE $g := \nu_i(\textbf{x})$
            \IF{$g>0$}
                \STATE $Q[i] := [a_i, c]$
            \ELSE
                \STATE $Q[i] := [c, b_i]$
            \ENDIF
        \ENDFOR
	    \STATE $\textbf{x} := (\frac{a_1+b_1}{2} \;\; \dotsi \;\; \frac{a_n+b_n}{2})^\top$
    \ENDWHILE
    \RETURN $\textbf{x}$
    \end{algorithmic}
\end{algorithm}

Обсудим корректность предложенного алгоритма. Далее, будем использовать следующие обозначения. Пусть на текущем уровне рекурсии рассматривается множество $Q\subset \mathbb{R}^n$ и некоторое его сечение $Q_k$ вида $Q_k = \left\{x \in Q\Big| x_k=c\right\}$ для некоторого $c\in \mathbb{R}$ и $k=\overline{1,n}$. Пусть $\nu$ есть вектор в $\mathbb{R}^n$. Тогда определим $\nu_{\parallel Q_k}$, $\nu_{\perp Q_k}$~--- проекции вектора $\nu$ на $Q_k$ и его ортогональное дополнение в пространстве $\mathbb{R}^{n}$ соответственно. Заметим, что $\nu_{\perp Q_k}$ есть скаляр. 

Сформулируем следующий вспомогательный результат.

\begin{lemma}\label{subgradient}
Пусть $f$~--- непрерывно-дифференцируемая выпуклая функция, и рассматривается задача её минимизации на множестве $Q_k \subset Q$. Если $\textbf{x}_*$~--- решение этой задачи, то существует условный субградиент функции $f$ на множестве $Q_x$ $g \in \partial_Q f(\textbf{x}_*)$
, такой что $g_\parallel = 0.$
\end{lemma}

Отметим, что такой Алгоритм \ref{alg:Dichotomy} сходится не для всех выпуклых функций, даже если допустить, что все вспомогательные подзадачи одномерной минимизации решаются точно. Отметим в этой связи пример негладкой выпуклой функции из \cite{Ston_Pas}, для которой теряется сходимость многомерной дихотомии.

Следующий результат~--- обобщение на многомерный случай утверждений \cite{Ston_Pas} о сходимости метода Ю.\,Е.\,Нестерова на квадрате.

\begin{theorem}\label{InexGradConst}
Пусть функция $f$ выпукла, имеет $L_{f}$-липшицев градиент при некотором $L_{f} > 0$. Пусть $\nu(\textbf{x}) = \nabla f(\textbf{x})$ и $\nu_{\perp Q_k}(\textbf{x}) = \left(\nu(\textbf{x})\right)_{\perp Q_k}$ для всякой текущей точки $\mathbf{x}$. Если все вспомогательные подзадачи решаются со следующей точностью по функции
\begin{equation}
\label{ConstEstFunc}
    \widetilde{\varepsilon} \leqslant \frac{\mu_f \varepsilon^2}{128 L_f^2 R^2},
\end{equation}
то после каждого удаления некоторой части (согласно пп. 11-14 Алгоритма \ref{alg:Dichotomy}) допустимого множества оставшаяся его часть содержит решение исходной задачи $\textbf{x}_*$ на гиперкубе $Q_x$.
\end{theorem}

Данная оценка требует точности решения вспомогательной задачи по функции порядка $\varepsilon^{2k}$ на $k$-ом уровне рекурсии. Таким образом, с учётом максимальной глубины рекурсии ($n-1$) получаем, что в худшем случае на каждом шаге алгоритма потребуется решать задачу с точностью $\varepsilon^{2n-2}$ по функции.

\begin{remark}\label{RemStrongConvDihot}
Заметим, что условие $\mu_f$-сильной выпуклости $f$ нужно лишь для теоретической оценки достаточной точности решения вспомогательных подзадач на итерациях алгоритма~\ref{alg:Dichotomy}, которая бы позволила гарантировать достижение требуемого качества решения задачи \eqref{func_f_1} за линейное время. Для практической реализации метода многомерной дихотомии ни знать константу $\mu_f$, ни требовать $\mu_f > 0$ нет необходимости, что существенно  использовано при постановке экспериментов.
\end{remark}

Интуитивно понятно, что для функции с липшицевым градиентом <<большая>> величина ортогональной компоненты при близости к точному решению (вспомогательной подзадачи) не сможет сильно уменьшиться, а следовательно, и изменить направление. С другой стороны, если эта компонента мала, и решается вспомогательная задача на многомерном параллелепипеде $Q_k$, то возможно выбрать эту точку в качестве искомого решения задачи на $Q_k$. Предложим некоторый подход к выбору точности для вспомогательных подзадач на базе отмеченной выше идеи.

Начнем с результата о необходимой точности для решения вспомогательных подзадач, гарантирующей сохранение искомого точного решения в допустимой области после удаления её частей на итерациях метода. Будем обозначать далее через $\Delta$ точность решения вспомогательных подзадач (п. 9 Алгоритма \ref{alg:Dichotomy}) по аргументу (п.\;10 алгоритма~\ref{alg:Dichotomy}).

\begin{theorem}\label{CurGrad}
Пусть функция $f$ выпукла и имеет $L_f$-липшицев градиент при $L_f > 0$, $\textbf{x}$~--- приближённое решение вспомогательной подзадачи на некоторой итерации метода.
Пусть на каждой итерации $\nu(\textbf{x}) = \nabla f(\textbf{x})$ и $\nu_{\perp Q_k}(\textbf{x}) = \left(\nu(\textbf{x})\right)_\perp$, а также приближение $\textbf{x}$ удовлетворяет условию:
$$\Delta \leqslant \frac{|\nu_{\perp Q_k}(\textbf{x})|}{L_f},$$
тогда после каждого удаления части допустимого множества оставшаяся его часть содержит решение исходной задачи $\textbf{x}_*$ на гиперкубе.
\end{theorem}

Полученная в теореме~\ref{CurGrad} оценка может привести к крайне медленной скорости сходимости Алгоритм \ref{alg:Dichotomy}, если проекция вектора-градиента во всякой текущей точке на ортогональное дополнение к рассматриваемому подпространству стремительно убывает при приближении к точке-решению. Поэтому сформулируем результат с альтернативным условием останова для вспомогательных задач. Обозначим через $Q_k$ некоторое подмножество (отрезок, часть плоскости или гиперплоскости) исходного гиперкуба $Q$, на котором решается вспомогательная подзадача.

\begin{theorem}\label{small}
Пусть функция $f$ выпукла, $M_f$-лип\-ши\-це\-ва и имеет $L_f$-лип\-ши\-цев градиент ($M_f, L_f > 0$). Тогда для достижения точности $\varepsilon$ по функции решения задачи \eqref{func_f_1} на множестве $Q_x$ достаточно выполнения следующего условия для всякого приближённого решения $\textbf{x}\in Q_k\subset Q$:
$$\Delta \leqslant \frac{\varepsilon- R|\nu_{\perp Q_k}(\textbf{x})|}{L_f+M_f R}, $$
где $R=a\sqrt{n}$~--- длина диагонали в исходном гиперкубе $Q_x$.
\end{theorem}

Здесь под $\varepsilon$ подразумевается точность по функции, с которой нужно решить задачу оптимизации на $Q_k$ при условии, что алгоритм находится на $k$-ом уровне рекурсии, то есть решает вспомогательную задачу на множестве $Q_k\subset Q$.

Объединяя полученные оценки, получаем следующую теорему для случая $n=2$.

\begin{theorem}\label{FullCond}
Пусть функция $f$ выпукла, $M_f$-лип\-ши\-це\-ва и имеет $L_f$-лип\-ши\-цев градиент ($M_f > 0, L_f >0$), а также в произвольной текущей точке $x$ при реализации метода мы используем неточный градиент $\nu(x)$, удовлетворяющий условию
\begin{equation}\label{tilde_delta}
\|\nabla f(\textbf{x}) - \nu(\textbf{x})\|_2 \leqslant \widetilde{\delta}(\textbf{x}).
\end{equation}
Тогда для того, чтобы после каждого удаления части допустимого множества оставшаяся его часть содержала решение исходной задачи $\textbf{x}_*$ на множестве $Q$, достаточно выполнения следующего условия для решения $\textbf{x}$ вспомогательной задачи минимизации $f$ на $Q_k$:
\begin{equation}\label{Adapt_inexact}
C_f\widetilde{\delta}(\textbf{x}) + \Delta \leqslant \max\left\{
	\frac{|\nu_{\perp Q_k}(\textbf{x})|}{L_f},
	\frac{\varepsilon - R |\nu_{\perp Q_k}(\textbf{x})|}{M_f+L_f R}
	\right\},
\end{equation}
где $C_f = \max\left(\frac{1}{L_f}, \frac{R}{M_f+L_f R}\right).$ В таком случае приближение искомого минимума с точностью $\varepsilon$ будет достигнуто не более чем за
\begin{equation}\label{dich_x_N}
N^* := \left\lceil\log_2 \left(\frac{4R(M_f+2L_f R)}{L_f\varepsilon}\right)\right\rceil
\end{equation}
итераций алгоритма~\ref{alg:Dichotomy}.
\end{theorem}

Отметим, что критерий останова \eqref{Adapt_inexact} применим ко внутренним подзадачам. Для внешних подзадач остановка определяется требуемым количеством итераций \eqref{dich_x_N}.

Кроме этого, заметим, что с учётом требований малости диаметра оставшейся части допустимой области после удаления гиперпрямоугольников, согласно предложенному алгоритму, оценка количества итераций, гарантирующих достижение $\varepsilon$-точ\-ности по функции, принимает следующий вид:

\begin{equation}\label{estimate_gen_dih}
O\left(\log_2\left(\varepsilon^{-1}\max\left(M_f, \frac{4(M_f +2L_f R)}{L_f}\right)\right)\right).
\end{equation}

Таким образом, на каждом уровне рекурсии (для гиперпрямоугольника $Q$) вспомогательная задача решается до тех пор, пока не будет выполнено условие \eqref{Adapt_inexact}. Причём в случае, если для некоторой точки условие \eqref{Adapt_inexact} выполняется и для самого второго аргумента максимума, то эта точка будет решением задачи на $Q$. При этом, если задача на $Q$ есть вспомогательная задача для некоторого гиперпрямоугольника $Q_1$ в более высокой размерности, то на более высоких уровнях рекурсии решение задачи, тем не менее, не останавливается.

Переходим теперь к описанию теоретических результатов об оценках скорости сходимости предложенной методики для сильно выпукло-вогнутых седловых задач вида \eqref{max_min} с достаточно малой размерностью (до 5) одной из групп переменных. Напомним, что мы рассматриваем исходную задачу \eqref{max_min} как задачу минимизации вспомогательного выпуклого функционала max-типа с использованием неточного градиента на каждой итерации (точность его нахождения регулируется за счёт решения вспомогательных минимизационных подзадач по другой группе переменных). Отметим, что поскольку множества $Q_x$ и $Q_y$ компактны, то $f$ удовлетворяет условию Липшица в силу предположений \eqref{smooth_F_11}--\eqref{smooth_F_44}.

\subsubsection{Оценки сложности алгоритма для седловых задач в случае использования многомерной дихотомии для подзадач малой размерности} 

Для минимизации по переменной малой размерности $x$ в задаче \eqref{max_min} будем использовать метод многомерной дихотомии~\ref{alg:Dichotomy} с неточным (аддитивно зашумленным) градиентом. Опишем необходимые условия на точность градиента целевого функционала, который будет использоваться на каждой итерации. Заметим, что согласно теореме Демьянова--Данскина (см. \cite{DDR-theorem}) имеем
$$
    \nabla f(\textbf{x}) = \nabla_x S(\textbf{x}, \textbf{y}(\textbf{x})).
$$
Здесь и далее будем считать, что $\nu(\textbf{x})=\nabla_x S(\textbf{x}, \textbf{y}_{\widetilde{\delta}}),$ где $\textbf{y}_{\widetilde{\delta}}$ есть приближение $\textbf{y}(\textbf{x})$ такое, что выполнено условие \eqref{tilde_delta}. Рассмотрим два возможных условия для вычисления $\textbf{y}_{\widetilde{\delta}}$.

1. В силу сделанных выше предположений для $S$:
$$
    \|\nu (\textbf{x})- \nabla f(\textbf{x})\|_2\leqslant L_{yx} \|\textbf{y}(\textbf{x}) - \textbf{y}_\delta\|_2.
$$

Таким образом, для того чтобы правильно определить оставшееся множество, достаточно выполнения неравенства \eqref{Adapt_inexact}, где $\widetilde{\delta}(\textbf{x}) = L_{yx} \|\textbf{y}(\textbf{x}) - \textbf{y}_\delta\|_2$.

Указанный выше способ предполагает возможность решения с линейной скоростью внутренней подзадачи задачи с любой точностью по аргументу. Однако если $S(\textbf{x}, \textbf{y})$~--- $\mu_y$-сильно вогнута по $\textbf{y}$, то эти условия можно заменить на условия сходимости по функции. В таком случае получаем, что вспомогательную задачу по $y$ нужно решать с точностью по функции $\delta$, удовлетворяющей условию
\begin{equation}\label{first_delta}
\delta \leqslant \frac{\mu_y}{2L_{yx}^2} \widetilde{\delta}^2.
\end{equation}

2. Можно получить другую оценку для необходимой точности решения вспомогательных подзадач. Заметим, что если $y_{\widetilde{\delta}}$ есть решение задачи максимизации вида $$f(x)=\max_{\textbf{y}\in Q_y}S(\textbf{x}, \textbf{y})$$
при фиксированном $x$ с точностью ${\delta}$ по функции, то $\nu(\textbf{x})$ есть ${\delta}$-градиент функции $f$ в точке $\textbf{x}$.
Тогда по теореме~\ref{lem:boundary}, в случае если расстояние от текущей точки $x$ до границы допустимого множества $Q_x$ достаточно велико, то есть $\rho(\textbf{x}, \partial Q)\geqslant 2\sqrt{\frac{\delta}{L_f}}$, верно следующее неравенство:

\begin{equation}
\|\nu(\textbf{x}) - \nabla f(\textbf{x})\|_2\leqslant 2\sqrt{{L_f\delta}}.
\end{equation}

В таком случае, для того чтобы правильно определить оставшееся множество, достаточно выполнения неравенства \eqref{Adapt_inexact} для $\widetilde{\delta}(\textbf{x}) = 2\sqrt{{L_f\delta}}$, где $\delta$~--- точность по функции решения вспомогательной подзадачи $y_{\widetilde{\delta}}$.

В данном случае получаем, что вспомогательную задачу по $y$ нужно решать с точностью по функции, удовлетворяющей условию
\begin{equation}
\label{second_delta}
\delta \leqslant \frac{2}{L_{f}} \widetilde{\delta}^2.
\end{equation}

Таким образом, на каждом шаге многомерной дихотомии будем вычислять неточный градиент согласно условиям \eqref{first_delta} или \eqref{second_delta}.

Далее, рассмотрим вспомогательную внутреннюю задачу минимизации по многомерной переменной. Стратегия решения вспомогательной задачи минимизации по многомерной переменной описана в разделе~\ref{second_x}. Точнее говоря, для вспомогательных подзадач максимизации по $y$ большой размерности будем использовать следующие методы в зависимости от выделенного класса задач (условий на $h$ или $S$).
\begin{enumerate}
    \item В случае проксимально-дружественной (\textbf{случай 1}) $h$ применяем быстрый градиентный метод для задач композитной оптимизации (алгоритм~4).
    \item Если $h$ имеет $L_h$-липшицев градиент (\textbf{случай 2}), то возможно применять ускоренный метод с реализацией разделения оракульных сложностей (алгоритм~2).
\end{enumerate}

В таких случаях получаем следующие оценки достаточного для достижения $\varepsilon$-точного решения задачи \eqref{max_min}, согласно определению~\ref{definepsaccuracy}, количества обращений к соответствующим вспомогательным подзадачам:

\begin{enumerate}
    \item В случае 1~---
    
        $O \left( 2^{n^2} \sqrt{\frac{L_{yy}}{\mu_y}} \log_2^{n+1} \left(\frac{CR}{\varepsilon}\right) \right)$ вычислений  $\nabla_y F(x, y)$  и решений подзадачи \eqref{prox_h}.
    \item В случае 2~---
    \begin{equation}\label{N_h}
    \begin{aligned}
        & O \left(  \sqrt{\frac{L_h}{\mu_y}} \log_2^{n+1}\left(\frac{CR}{\varepsilon}\right) \right) \text{вычислений } \nabla h(y) \\
        & \text{и } O \left( 2^{n^2} \left( \sqrt{\frac{L_h}{\mu_y}}+\sqrt{\frac{L_{yy}}{\mu_y}} \right) \log_2^{n+1}\left(\frac{CR}{\varepsilon}\right) \right) \text{ вычислений } \nabla_y F(x, y).
    \end{aligned}
    \end{equation}
\end{enumerate}

Отметим, что условиями останова для решения вспомогательной задачи являются условия \eqref{first_delta} и \eqref{second_delta}.

Также заметим, что в случае сепарабельной функции аналогичным образом, как и ранее, в пункте~\ref{second_x}, можно гарантировать достижение $\varepsilon$-решения задачи \eqref{max_min} за $O \left( n \ln \frac{n}{\varepsilon} \right)$ вычислений $\nabla_x F$ и $O \left( mn \ln \frac{n}{\varepsilon} \ln \frac{m}{\varepsilon} \right)$ вычислений $S(x, y)$ при ослабленных условиях на гладкость.

\section{Результаты вычислительных экспериментов}\label{labsect3}

\subsection{Постановка задач, для которых проводится сравнение вычислительной эффективности предложенных методов}

В качестве важного класса седловых задач можно выделить лагранжевы седловые задачи, связанные с задачами выпуклого программирования. Если у такой задачи два функционала ограничения и выполнено условие Слейтера, то двойственная задача максимизации двумерна и к ней вполне возможно применить метод многомерной дихотомии (двумерный случай) после локализации допустимой области значений двойственных переменных. Более того, полученные нами результаты теоретически обосновывают линейную скорость сходимости для случая гладкого сильно выпуклого функционала и выпуклого гладкого функционала ограничения. Похожее верно, если также, например, применять метод эллипсоидов для двойственной задачи. Однако возможна ситуация, когда за счёт отсутствия необходимости на итерациях находить точное значение градиента целевой функции метод Ю.\,Е.\,Нестерова \cite{Ston_Pas} работает быстрее, причём даже в случае негладких функционалов ограничений max-типа.

Рассмотрим задачу вида
\begin{equation}\label{eq_sedlo_problem}
	\max_\lambda \left\{ \phi(\lambda) = \min_x F(x, \lambda)\right\}.
\end{equation}
Оптимальную точку $\lambda_*$ будем находить, используя метод многомерной дихотомии, решая на каждом шаге вспомогательную задачу минимизации по $x$ быстрым градиентным методом. Подробно про это сказано в разделе~\ref{theor_dich}.

Ниже выделим несколько исследованных в настоящей работе подходов, с использованием которых можно решить задачу \eqref{eq_sedlo_problem}. 

Во-первых, можно решать основную задачу при помощи метода эллипсоидов при решении вспомогательной задачи быстрым градиентным методом (см. разд.~\ref{1st_approach}, случай 2).

Во-вторых, возможно решить поставленную задачу при помощи быстрого градиентного метода с $(\delta, L)$-оракулом \cite{FGM,th2_cite,devolder_for_dLm_model}, при этом решая вспомогательную задачу при помощи метода эллипсоидов (алгоритм~\ref{alg:ellipsoid}). 

В-третьих, к рассматриваемой задаче можно подойти с использованием метода многомерной дихотомии, решая вспомогательную задачу быстрым градиентным методом (см. разд.~\ref{3st_approach}, случай 2).

Наконец, в четвёртом способе методика не будет учитывать малую размерность одной из переменных, и будет использоваться вариант быстрого градиентного метода с $(\delta,L)$-оракулом \cite{FGM,th2_cite,devolder_for_dLm_model} целевой функции для внешней задачи и обычный быстрый градиентный метод (БГМ) для внутренней подзадачи.

Стоит заметить, что в предыдущих подразделах мы рассматривали сильно выпукло-вогнутые седловые задачи вида \eqref{problem:min_max00}. Тем не менее, сильная выпуклость оптимизационных подзадач небольшой размерности (к которым мы применяем методы секущей гиперплоскости или авторский вариант многомерной дихотомии) существенна лишь для вывода теоретических оценок сложности. Практическая же реализация алгоритмов~\ref{alg:ellipsoid},~\ref{alg:vaidya} и~\ref{alg:Dichotomy} возможна и без предположения о сильной выпуклости. В проведённых нами экспериментах удаётся без этих предположений подобрать вспомогательные параметры указанных методов, без использования сильной выпуклости задачи \eqref{problem:min_max00} по переменным небольшой размерности (в нашем случае это двойственные переменные лагранжевых седловых задач), чтобы достигнуть желаемого качества приближённого решения. Этим объясняется корректность приведённых в данном разделе результатов экспериментов, несмотря на то что рассматриваемые задачи не сильно выпуклы (вогнуты) по одной из групп переменных (а только выпуклы или вогнуты).


\subsection{Лагранжева седловая задача к задаче квадратичной оптимизации}\label{section:experiments_triangle}

В качестве конкретного примера для сравнительных расчётов рассмотрим задачу квадратичной оптимизации с ограничениями:

\begin{equation}
\label{quadratic_problem}
\min_{\substack{x \in Q_r \subset \mathbb{R}^n \\ g^i(x) \leqslant 0, \; i = 1,\ldots,m}} \left\{ f(x) := \frac{1}{2} \|Ax - b\|_2^2 \right\},
\end{equation}
где $A \in \mathbb{R}^{n \times n}$, $b \in \mathbb{R}^n$, и $Q_r = \left\{x\;\Big|\;\|x\|_2 \leqslant r\right\}$ --– евклидов шар, а каждое из ограничений $g^i(x)$ является линейным:

$$
    g^i(x) = {c^i}^\top x + d^i,\quad c^i \in \mathbb{R}^n, d^i \in \mathbb{R}.
$$

В дальнейшем для возможности использования метода оптимизации на квадрате или треугольнике мы будем работать с двумя выпуклыми негладкими ограничениями $g_1$ и $g_2$, $\max$-агрегирующими исходные ограничения:

$$
g_1(x) = \max \left\{ g^i(x)\;\Big|\;i = 1, \ldots, \left\lfloor\frac{m}{2}\right\rfloor \right\},
$$
$$
g_2(x) = \max \left\{ g^i(x)\;\Big|\;i = \left\lfloor\frac{m}{2}\right\rfloor + 1, \ldots, m\right\}.
$$

В такой постановке с двумя ограничениями исходная задача \eqref{quadratic_problem} будет иметь двойственную задачу следующего вида:

$$\max_{\lambda_1 + \lambda_2 \leqslant \Omega_\lambda} \left\{ \varphi(\lambda_1, \lambda_2) := \min_{x \in Q_r} \left\{f(x) + \lambda_1 g_1(x) + \lambda_2 g_2(x) \right\} \right\},$$
где константа $\Omega_\lambda$ оценивается из условия Слейтера следующим образом:
$$\Omega_\lambda = \frac{1}{\gamma} f(\hat{x}),$$
где $\gamma = - \max\left\{g_1(\hat{x}), g_2(\hat{x})\right\} > 0$, а $\hat{x}$ --– внутренняя точка множества, задаваемого исходными ограничениями. Таким образом, вместе с условием неотрицательности двойственных переменных мы получаем, что множество, на котором решается двойственная задача, представляет из себя прямоугольный треугольник с катетами длины $\Omega_\lambda$, лежащими на осях координат.

Пусть $A$ –- разреженная матрица с долей ненулевых элементов $\sigma$ со случайными положительными диагональными элементами из равномерного распределения $A_{ii} \propto \mathcal{U}(0, 1.1)$ и ненулевыми недиагональными элементами также из равномерного распределения $A_{ij} \neq 0, A_{ij} \propto \mathcal{U}(0, 1)$, элементы вектора $b$~--- независимые равномерно распределенные случайные величины $b_i \propto \mathcal{U}(0, 0.5)$, вектор $c^i$ и скаляр $d^i$, задающие $i$-е ограничение, также генерируются случайно из равномерного распределения $\mathcal{U}(0, 0.1)$.

Сравним скорость работы метода двумерной дихотомии (далее – метод оптимизации на квадрате), метода оптимизации на (равнобедренном прямоугольном) треугольнике (который аналогичен методу двумерной дихотомии)  на множестве $Q = \left\{x \in \mathbb{R}_{++}^2 \;\Big|\; x_1 + x_2 \leqslant \Omega_\lambda\right\}$, метода эллипсоидов и быстрого градиентного метода. Приведём здесь описание используемого далее варианта метода дихотомии на равнобедренном прямоугольном треугольнике. Каждая его итерация осуществляется в соответствии со следующим алгоритмом.
\begin{enumerate}
    \item На первом шаге итерации метода проводится разделяющий отрезок, соединяющий середины одного из катетов и гипотенузы (рис.~\ref{fig:triangle_method_1}). С помощью одного из методов одномерной оптимизации (например, метода золотого сечения) решается вспомогательная задача минимизации на данном отрезке.
    
    \item Далее, в полученной точке решения вспомогательной задачи $x$ вычисляется градиент $\nabla f(x)$, после чего отсекается та часть треугольника, в которую он направлен.
    
    \item В случае, если после отсечения остаётся в два раза меньший гомотетичный исходному треугольник, метод переходит к следующей итерации.
    
    \item Иначе, на втором шаге итерации аналогичным образом отсекается одна из частей оставшейся трапеции, на которые её делит отрезок, соединяющий середины гипотенузы и другого катета исходного треугольника (рис.~\ref{fig:triangle_method_2}). Если в результате этого остаётся треугольник, то метод переходит к следующей итерации. Если же остаётся квадрат, то для дальнейшей оптимизации используется метод двумерной дихотомии.
\end{enumerate}

\begin{figure}
     \centering
     \begin{subfigure}[b]{0.38\textwidth}
         \centering
         \includegraphics[width=0.7\textwidth]{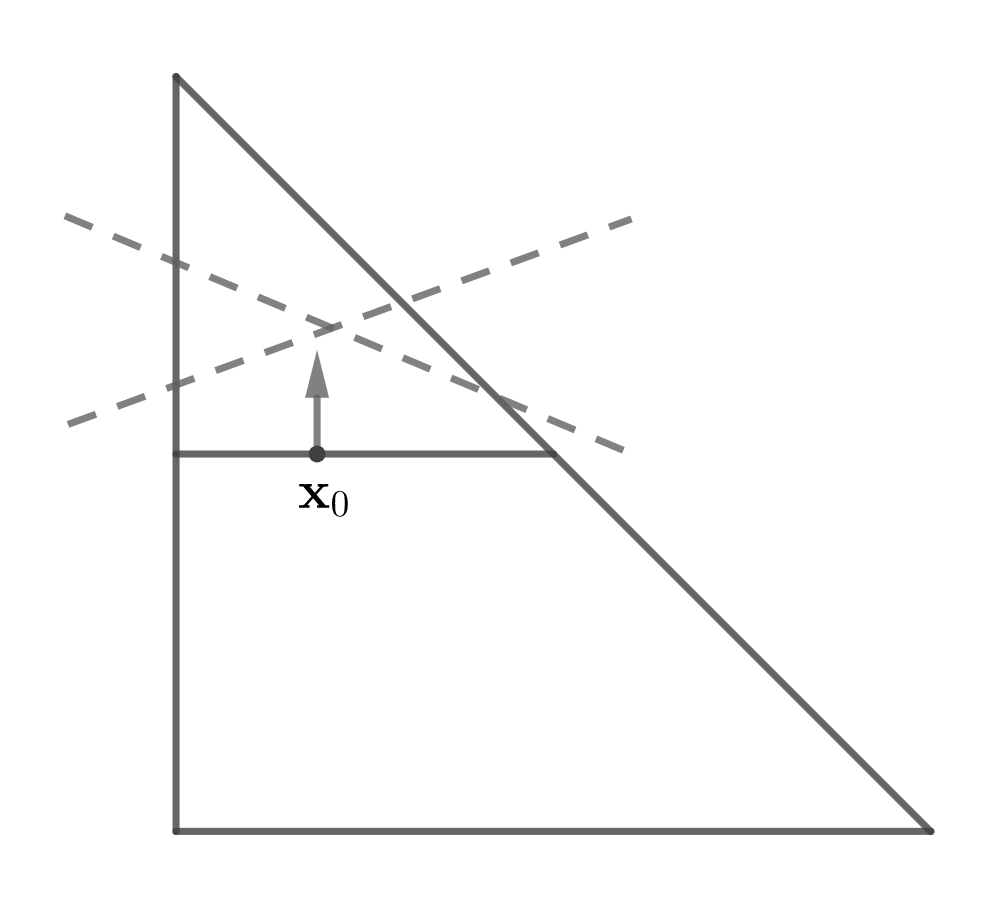}
         \caption{Первый шаг отсечения}
         \label{fig:triangle_method_1}
     \end{subfigure}
     \hfill
     \begin{subfigure}[b]{0.6\textwidth}
         \centering
         \includegraphics[width=0.68\textwidth]{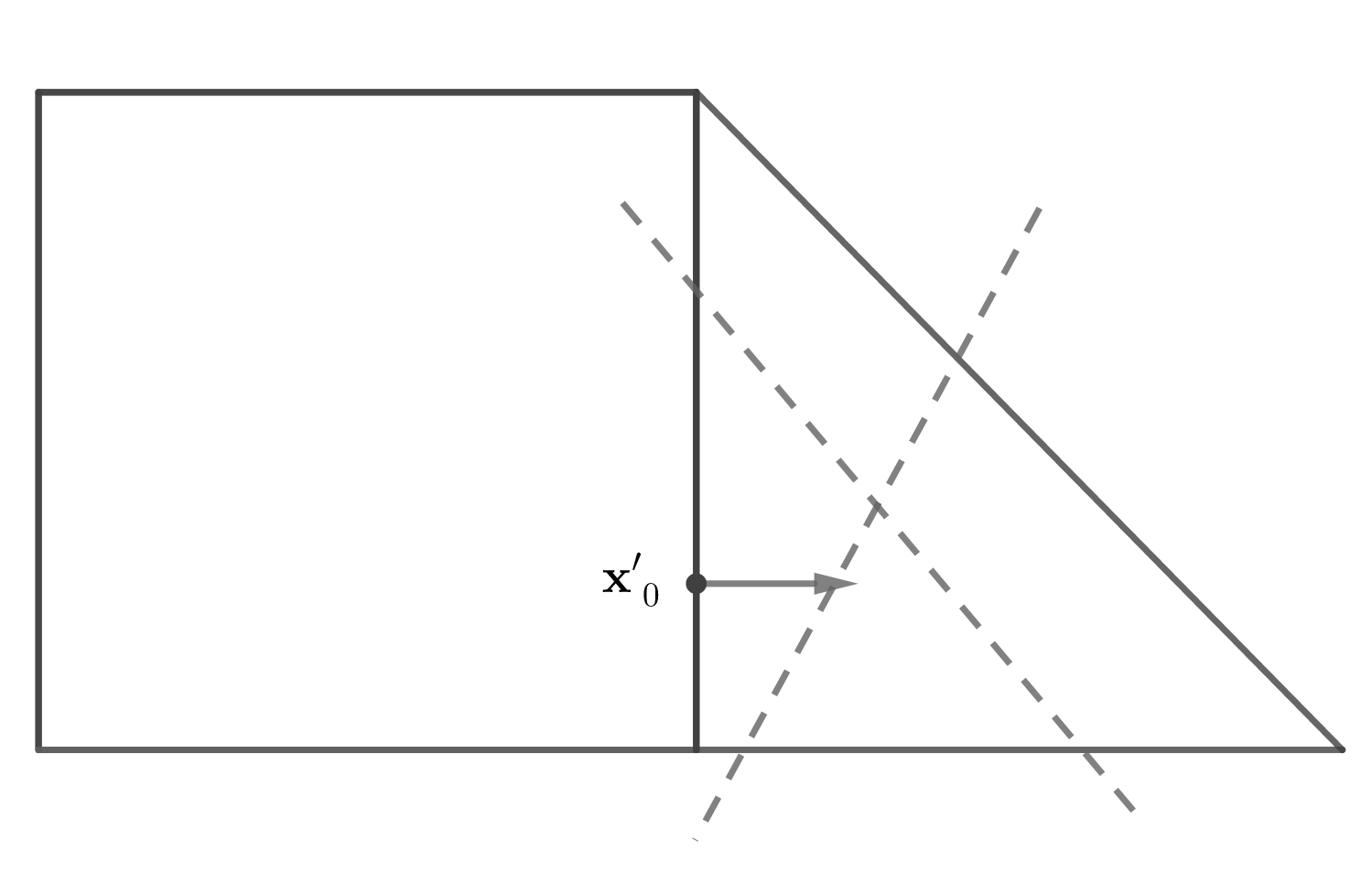}
         \caption{Второй шаг отсечения}
         \label{fig:triangle_method_2}
     \end{subfigure}
        \caption{Иллюстрация к описанию метода дихотомии на треугольнике.}
        \label{fig:triangle_method}
\end{figure}

В качестве критерия останова для всех сравниваемых в настоящем пункте методов будем использовать условие:

$$|\lambda_1 g_1(x_\delta(\lambda)) + \lambda_2 g_2(x_\delta(\lambda))| < \varepsilon,$$
где $x_\delta(\lambda)$~--- приближенное с точностью $\delta$ по функции решение вспомогательной задачи $\min_{x \in Q} \left\{f(x) + \lambda_1 g_1(x) + \lambda_2 g_2(x) \right\}$. Выполнение данного условия гарантирует достижение точности по функции решения исходной задачи

$$f(x_\delta(\lambda)) - \min_{\substack{x \in Q \\ g_1(x), g_2(x) \leqslant 0}} f(x) \leqslant \varepsilon + \delta.$$

При этом на каждой итерации метода двойственные множители положительны в силу выбранной области их локализации (треугольник в положительном ортанте), что исключает возможность преждевременного выполнения этого условия в точке 0.

Вспомогательная задача минимизации решается с помощью субградиентного метода \cite{Shor79}. Число итераций субградиентного метода определяется экспериментально таким образом, чтобы в точке решения задачи достигалась точность по функции $\delta = 0.005$ (при сравнении с решением, полученным на большом числе итераций метода), а также так, чтобы значения ограничений $g_1$ и $g_2$ в данных точках были неположительными. Для $n = 400$ при $r = 5$ и $\sigma = 0.005$ оказалось достаточно $800$ итераций субградиентного метода, а для $n = 1000$, $r = 2$ и $\sigma=0.001$~--- $2500$ итераций.

Как можно видеть из таблицы~\ref{tab1}, где для различных значений $\varepsilon$ сравнивается работа методов для $n=400$ при значениях $m=10,\,m=20$, предложенный метод оптимизации на треугольнике достигает выполнения критерия останова, а значит, и заданной точности по функции за меньшее по сравнению с методом эллипсоидов число итераций и время работы.

\begin{center}
\begin{table}[ht!]
\centering
\begin{tabular}{|c|c|c|c|c|c|}
\hline
\multirow{2}{*}{$\varepsilon$}& \multirow{2}{*}{\;\;\;$m$\;\;\;}
& \multicolumn{2}{c|}{Метод на треугольнике}
& \multicolumn{2}{c|}{Метод эллипсоидов} \\ \cline{3-6}
&   & Итерации & Время, мс & Итерации & Время, мс  \\ \hline
\multirow{2}{*}{0.5}  & 10  & 4     & 687  & 8      & 820  \\ \cline{2-6}
                      & 20  & 4     & 818  & 8      & 1170 \\ \hline
\multirow{2}{*}{0.1}  & 10  & 8     & 810  & 14     & 1300 \\ \cline{2-6}
                      & 20  & 8     & 1020 & 14     & 1390 \\ \hline
\multirow{2}{*}{0.05} & 10  & 12    & 1160 & 16     & 1460 \\ \cline{2-6}
                      & 20  & 14    & 1840 & 18     & 2140 \\ \hline
\multirow{2}{*}{0.01} & 10  & 22    & 3170 & 36     & 3530 \\ \cline{2-6}
                      & 20  & 26    & 3260 & 38     & 3770 \\ \hline
\end{tabular}
\caption{Сравнение работы методов при $n = 400$.}
\label{tab1}
\end{table}
\end{center}

Из таблицы~\ref{tab2} сравнения методов для $n = 1000$ и $m = 10$ также видно, что предложенная схема эффективнее метода эллипсоидов. При этом время работы метода оптимизации на треугольнике несколько меньше, чем время работы метода оптимизации на квадрате, что достигается благодаря локализации двойственных множителей на прямоугольном треугольнике и вследствие этого меньшей длине отрезков, на которых необходимо решать дополнительные одномерные оптимизационные задачи на первых итерациях методов. Кроме того, значительно больших в сравнении с прочими методами числа итераций и времени работы требует метод эллипсоидов для задачи с $m$ ограничениями (не агрегированными в два ограничения вида $\max$) ввиду увеличения размерности задачи и больших временных затрат на выполнение матрично-векторных операций.

\begin{center}
\begin{table}[ht!]
\centering
\scalebox{0.75}{
\begin{tabular}{|c|c|c|c|c|c|c|}
\hline
\multirow{2}{*}{$\varepsilon$} & \multicolumn{2}{c|}{Метод на квадрате} & \multicolumn{2}{c|}{Метод на треугольнике}
& \multicolumn{2}{c|}{Метод эллипсоидов}\\ \cline{2-7}
& Итерации & Время, с & Итерации & Время, с & Итерации & Время, с  \\ \hline
$0.5$   & 6  & 6.10 & 6 & 5.41  & 4   & 6.12  \\ \hline
$0.1$   & 12 & 8.92 & 12& 8.25  & 16  & 12.7  \\ \hline
$0.05$  & 18 & 12.6 & 16& 11.2  & 24  & 23.8   \\ \hline
$0.01$  & 24 & 25.3 & 22& 24.1  & 30  & 32.5   \\ \hline
\multirow{2}{*}{$\varepsilon$} & \multicolumn{2}{c|}{Метод эллипсоидов ($m=10$)}  & \multicolumn{2}{c|}{Метод эллипсоидов ($m=10$, БГМ)} \\ \cline{2-5}
        & Итерации & Время, с & Итерации & Время, с  \\ \cline{1-5}
$0.5$   & 8        & 15.3     & 6        & 6.27      \\ \cline{1-5}
$0.1$   & 20       & 26.2     & 10       & 17.6      \\ \cline{1-5}
$0.05$  & 32       & 38.7     & 34       & 38.8      \\ \cline{1-5}
$0.01$  & 40       & 49.5     & 40       & 50.2      \\ \cline{1-5}
\end{tabular} }
\caption{Сравнение работы методов при $n = 1000, m = 10$.}
\label{tab2}
\end{table}
\end{center}

Сравним работу метода на треугольнике для решения внешней задачи с быстрым градиентным методом при двух $\max$-агрегированных ограничениях и при $m=10$ исходных ограничениях (таблица~\ref{tab3}). В обоих вариантах быстрый градиентный метод требует для достижения той же точности большее число итераций и время работы, чем метод на треугольнике.

\begin{center}
\begin{table}[ht!]
\centering
\begin{tabular}{|c|c|c|c|c|}
\hline
\multirow{2}{*}{$\varepsilon$} & \multicolumn{2}{c|}{\;\;\;\;\;\;БГМ\;\;\;\;\;\;} & \multicolumn{2}{c|}{\;\;\;БГМ ($m=10$)\;\;\;}\\ \cline{2-5}
& Итерации & Время, с  & Итерации & Время, с  \\ \hline
$0.5$   & 10       & 10.8      & 12       & 12.3       \\ \hline
$0.1$   & 16       & 20.6      & 16       & 22.9       \\ \hline
$0.05$  & 22       & 34.1      & 22       & 34.8       \\ \hline
$0.01$  & 28       & 36.9      & 32       & 37.3       \\ \hline
\end{tabular}
\caption{Работа быстрого градиентного метода при $n = 1000, m = 10$.}
\label{tab3}
\end{table}
\end{center}

\subsection{Лагранжева седловая задача к задаче LogSumExp с линейными функционалами ограничений}

Рассмотрим LogSumExp задачу с $\ell_2$-регуляризацией и линейными ограничениями:

\begin{equation*}
\label{LogSumExp_problem}
\begin{split}
    \min_{x \in \mathbb{R}^m} \left\{\log_2 \left(1+\sum_{k=1}^m e^{\alpha x_k}\right) + \frac{\mu_x}{2} \| x \|_2^2 \right\}, \\
    \text{удовл.} \;\;Bx \leqslant c,\; B \in \mathbb{R}^{n \times m},\; c \in \mathbb{R}^{n},\; \alpha \in \mathbb{R}.
\end{split}
\end{equation*}

Лагранжиан этой задачи можно записать в следующей форме:
\begin{equation*}
    r(x) + F(x, y) - h(y),
\end{equation*}
где
\begin{equation*}
    r(x) = \frac{\mu_x}{2} \| x \|_2^2,\quad F(x, y) = \log_2 \left(1+\sum_{k=1}^m e^{\alpha x_k}\right) + y^T Bx,\quad h(y) = y^T c.
\end{equation*}
Тогда двойственная задача является выпукло-вогнутой седловой задачей:
\begin{equation}
\label{LogSumExpdual}
    \max_{y \in \mathbb{R}_+^n} \min_{x \in \mathbb{R}^m} \left\{r(x) + F(x, y) - h(y)\right\}.
\end{equation}

Заметим, что в приведенной постановке задачи $r(x)$ является проксималь\-но-дружественной.

Согласно теореме~\ref{constrains_y} (см. приложение) мы знаем, что
$$y_* \in Q_y = \left\{y \in \mathbb{R}_+^n\Big| y_k \leqslant \frac{f(x_0)}{\gamma}\right\},$$
где $x_0$ есть внутренняя точка полиэдра $Bx \leqslant c$, а $\gamma=\min_k\left\{c_k - [Bx]_k\right\} > 0.$ Также легко видеть, что $x$ должен лежать внутри шара $Q_x = B_{R_x}(0)$ с центром в нуле и некоторым конечным радиусом $R_x$. Действительно, значение функции в нуле есть $s_0 = S(0, y) = \log_2(m+1)-y^\top c$ для любого $y\in Q_y$, и можно найти такой $x$, что квадратичная часть по $x$ будет больше данного значения для любого $y \in Q_y$.

Обсудим параметры, связанные с константами Липшица градиентов рассматриваемых функций. Для функций $r$ и $h$ они очевидны:
$$L_r = \mu_x, \; L_h = 0.$$
Также очевидны следующие константы для функции $F$:
$$L_{x y} = \|B\|_2 R_y,\,L_{y x} = \|B\|_2 R_x,\,L_{y y} = 0.$$
Константа $L_{xx}$ есть сумма констант для LogSumExp, которую мы сейчас вычислим, и линейной функции по $x$, которая есть ноль. Константа Липшица для градиента LogSumExp есть максимальное собственное число ее гессиана, которое равно $\alpha$.  Таким образом,
$$L_{xx} = L_{LSE} = \max_{x} \lambda_1 \nabla^2 LSE(x) = \alpha,$$
где $LSE(x) = 
\log_2 \left(1+\sum_{k=1}^m e^{\alpha x_k}\right)$. 

Заметим, что в полученных в предыдущих пунктах теоретических оценках исследуемых методов требуется сильная вогнутость по $y$. В случае дифференцируемости внутренней функции, согласно теореме (\ref{dual_L_from_mu}), имеем сильную вогнутость с константой $\mu_y = \frac{M_g}{\mu_x},$ где $M_g = \|B\|_2$ есть константа Липшица для функции $g$ для случая, когда маломерная задача решается как основная. В случае, когда она решается как вспомогательная, используем прием регуляризации, то есть добавим квадратичный член $\frac{\varepsilon}{2R^2}$, где $R \geqslant \|y_*-y_0\|_2$, $y_0$~--- стартовая точка метода.

Параметры $\alpha_k$ генерируются из равномерного распределения $\mathcal{U} (-\alpha_0, \alpha_0)$, $\alpha_0 = 0.001$. Элементы матрицы $B$ генерируются из равномерного распределения $\mathcal{U} (-k, k),\,k = 10^3$. Параметр $\mu_x = 0.001$, элементы вектора $c$ равны единице.

Уточним использованное при проведении экспериментов условие останова сравниваемых методов. Заметим, что если $x_\delta(\lambda)$ есть решение с точностью $\delta$ по функции задачи
\begin{equation}\label{eqLagrangeSadle}
    \min_{x\in Q_x} \left\{ f(x) + \lambda^\top g(x) \right\},
\end{equation}
и при этом выполняется дополнительное условие
\begin{equation}
\label{constr_cond}
	\Big|\lambda^\top g(x_\delta(\lambda))\Big| \leqslant \varepsilon,
\end{equation}
то $x_\delta(\lambda)$~--- приближённое решение задачи минимизации $f$ с точностью $\delta + \varepsilon$ по функции, то есть
$$f(x_\delta(\lambda)) - \min_{x\in Q_x} f(x) \leqslant \delta + \varepsilon.$$
Действительно,
\begin{equation*}
\begin{aligned}
f(x_\delta(\lambda)) + \lambda^\top g(x_\delta(\lambda)) & \leqslant
f(x(\lambda)) + \lambda^\top g(x(\lambda)) + \delta \\& \leqslant 
\phi(\lambda_*) + \delta =
f(x_*) + \lambda_*^\top g(x_*) +\delta =f(x_*).
\end{aligned}
\end{equation*}
Последний переход совершён в силу условия Каруша--Куна--Таккера, утверждающего, что в точке $(\lambda_*, x(\lambda_*))$ должно быть выполнено условие дополняющей нежёсткости $\lambda_i g_i(x) = 0\; \forall i$.

Итак, возможно выбрать следующее условие останова:

\begin{equation}
\begin{cases}
	|\lambda^\top g(x_{\delta}(\lambda))| \leqslant \frac{\varepsilon}{2},\\
	g_i(x) \leqslant 0,\;\forall i:\;\lambda_i = 0.
\end{cases}
\end{equation}

Это условие гарантирует точность $\varepsilon$ по функции $f$, как уже выше было указано. Второе условие добавлено в силу того, что при $\lambda_i=0$ величина несоответствия условию $g_i(x)$ может быть сколь угодно большой.

Применяя условие \eqref{constr_cond} для останова метода при решении задачи \eqref{LogSumExpdual}, сравним описанные выше методы. Кроме этого, будет поставлено дополнительное ограничение на время выполнения метода. В случае, если предел времени превышен до того, как выполнено условие останова, то выполнение метода прерывается и возвращается текущий результат. В наших задачах мы установили этот предел равным $100$ секундам. Знак <<->> в соответствующих таблицах означает, что данный метод не успел завершиться за выделенное время при данных параметрах.

Для расчетов был использован Python версии 3.7.3 с устанавливаемой библиотекой numpy версии 1.18.3. Весь код выложен в репозиторий на GitHub (см. \cite{code_experiments}).

Результаты экспериментов для размерностей $n=2, 3, 4$ по двойственной переменной представлены в таблицах~\ref{res_m_2},~\ref{res_m_3} и~\ref{res_m_4}. В данных таблицах указано время работы для случая, когда маломерная задача решается быстрым градиентным методом или маломерными методами, такими как метод эллипсоидов с $\delta$-субградиентом или метод многомерной дихотомии, представленный в данной работе, а вспомогательная многомерная решается при помощи быстрого градиентного метода.

\begin{center}
\begin{table}[ht!]
\centering
\begin{tabular}{|c|c|c|c|c|c|}
\hline
\multirow{2}{*}{$\varepsilon$}& \multirow{2}{*}{\;\;\;$m$\;\;\;}
& \multicolumn{4}{c|}{Время работы, c}  \\ \cline{3-6}
& & БГМ & Эллипсоиды & Дихотомия & Вайда  \\ \hline
\multirow{2}{*}{$10^{-3}$}  & $10^2$  & \bf{0.02} & 0.27 & 0.14 & 0.39  \\
							& $10^3$  & \bf{0.03} & 0.53 & 0.27 & 0.56\\
                            & $10^4$  & \bf{0.45} & 9.86 & 4.33 & 6.98\\ \hline
\multirow{2}{*}{$10^{-6}$}  & $10^2$  & 3.48 & 0.45 &\bf{0.22} & 0.50\\
							& $10^3$  & 0.47 & 0.85  & \bf{0.45} & 0.79\\
                            & $10^4$  & \bf{0.63} & 16.72 & 6.16 & 11.10\\ \hline
\multirow{2}{*}{$10^{-9}$}  & $10^2$  & - & 0.79 & \bf{0.67} & 0.71\\
							& $10^3$  & - & 1.45 & \bf{1.12} & 1.24\\
                            & $10^4$  & - & 26.23 & \bf{16.09} & 16.82\\ \hline
\end{tabular}
\caption{Сравнение работы методов при $n = 2$.}
\label{res_m_2}
\end{table}
\end{center}

\begin{center}
\begin{table}[ht!]
\centering
\begin{tabular}{|c|c|c|c|c|c|}
\hline
\multirow{2}{*}{$\varepsilon$}& \multirow{2}{*}{\;\;\;$m$\;\;\;}
& \multicolumn{4}{c|}{Время работы, c}  \\ \cline{3-6}
& & БГМ & Эллипсоиды & Дихотомия & Вайда  \\ \hline
\multirow{2}{*}{$10^{-3}$}  & $10^2$  & \bf{0.05} & 0.65 & 0.88 & 0.71  \\
							& $10^3$  & \bf{0.03} & 1.30 & 1.56 & 0.91\\
                            & $10^4$  & \bf{0.36} & 22.05 & 20.52 & 10.92\\ \hline
\multirow{2}{*}{$10^{-6}$}  & $10^2$  & 2.46 & 1.07 &- & \bf{0.79}\\
							& $10^3$  & \bf{0.53} & 2.06 & - & 1.27\\
                            & $10^4$  & \bf{0.61} & 37.64 & - & 17.66\\ \hline
\multirow{2}{*}{$10^{-9}$}  & $10^2$  & - & 1.89 & - & \bf{1.17}\\
							& $10^3$  & - & 3.63 & - & \bf{1.64}\\
                            & $10^4$  & - & 59.71 & - & \bf{25.41}\\ \hline
\end{tabular}
\caption{Сравнение работы методов при $n = 3$.}
\label{res_m_3}
\end{table}
\end{center}

\begin{center}
\begin{table}[ht!]
\centering
\begin{tabular}{|c|c|c|c|c|c|}
\hline
\multirow{2}{*}{$\varepsilon$}& \multirow{2}{*}{\;\;\;$m$\;\;\;}
& \multicolumn{4}{c|}{Время работы, c}  \\ \cline{3-6}
& & БГМ & Эллипсоиды & Дихотомия & Вайда  \\ \hline
\multirow{2}{*}{$10^{-3}$}  & $10^2$  & \bf{0.06} & 1.08 & 7.06 & 0.9  \\
							& $10^3$  & \bf{0.03} & 2.22 & 12.24 & 1.38\\
                            & $10^4$  & \bf{0.37} & 40.37 & - & 16.06\\ \hline
\multirow{2}{*}{$10^{-6}$}  & $10^2$  & 3.10 & 1.86 & - & \bf{1.15}\\
							& $10^3$  & \bf{0.56} & 3.82 & - & 1.90\\
                            & $10^4$  & \bf{0.67} & - & - & 25.03\\ \hline
\multirow{2}{*}{$10^{-9}$}  & $10^2$  & - & 3.42 & - & \bf{2.01}\\
							& $10^3$  & - & 6.20 & - & \bf{2.83}\\
                            & $10^4$  & - & - & - & \bf{33.12}\\ \hline
\end{tabular}
\caption{Сравнение работы методов при $n = 4$.}
\label{res_m_4}
\end{table}
\end{center}

Случай, когда маломерная задача является вспомогательной, не представлен в таблицах, поскольку во всех экспериментах данные методы не успевали сойтись в рамках выделенного времени. Из этого можно сделать вывод, что при малых размерностях ограничений выгоднее решать маломерную задачу как основную.

Обсудим полученные результаты. Во-первых, для всех $n = 2, 3, 4$ (таблицы~\ref{res_m_2},~\ref{res_m_3} и~\ref{res_m_4} соответственно) видно, что при не очень высокой требуемой точности $\varepsilon = 10^{-3}$ быстрый градиентный метод показывает наилучший результат. Действительно, при данном значении точности он быстрее как минимум на порядок в сравнении с маломерными методами (методом эллипсоидов и методом дихотомии). С другой стороны, при повышении требуемой точности маломерные методы становятся быстрее градиентного метода. Так, например, для $n=2$ (см. таблицу~\ref{res_m_2}) мы видим, что для точности $\varepsilon = 10^{-9}$ маломерные методы быстрее многомерных для всех размерностей исходной задачи $m$.

Во-вторых, заметим, что при $n=2$ предложенный метод многомерной дихотомии сходится быстрее, чем метод эллипсоидов и метод Вайды для всех $\varepsilon$ и $m$. Однако сложность этого метода очень быстро растет с размерностью (см. оценку сложности \eqref{dich_x}), что проявляется и в экспериментах. Так, уже для $n=3$, то есть при повышении размерности на один он перестает быть эффективным по сравнению с другими методами. При $m>2$ и $\varepsilon=10^{-9}$ в большинстве тестов метод Вайды показывает себя лучше, чем остальные методы.

В-третьих, отметим характер зависимости времени работы от размерности прямой задачи $m$. Время работы при повышении $m$ и при фиксированных $n$ и $\varepsilon$ растет быстрее для маломерных методов, чем для многомерных. Следствием этого является то, что рассмотренные в работе методы маломерной минимизации работают эффективнее по сравнению с быстрым градиентным методом в наших экспериментах при $\varepsilon=10^{-6}$ только при малой размерности $m=100$.

\subsection{Лагранжева седловая задача к задаче LogSumExp с линейными функционалами ограничений при аддитивно зашумленном градиенте}\label{sec:exp-additive-noise}

В условиях предыдущего примера рассмотрим аналогичную постановку задачи \eqref{LogSumExpdual} и применим для её решения метод эллипсоидов для внешней $\max$-задачи вместе с быстрым градиентным методом для внутренней $\min$-задачи при дополнительном условии, что получение градиента по переменной $y$ внешней задачи при вызове соответствующего оракула в алгоритме~\ref{alg:ellipsoid} происходит с некоторой аддитивной ошибкой. Точнее говоря, если точное значение градиента функции под оператором $\min$ есть $\nabla_y f(y)$, вместо него доступен вектор $v$, такой, что:
\[
    \|\nabla_y f(y) - v\|_2 \leq \Delta.
\]
В соответствии с замечанием~\ref{inexactnesses_ellipsoid} $\Delta$-аддитивная неточность градиента может учитываться как дополнительная $\delta$-неточность оракула двумя способами: равномерно, используя $\mu$-сильную выпуклость (имеющую место для рассматриваемой задачи), или же динамически изменяя $\delta$ в зависимости от диаметра текущего эллипсоида, равного (в обозначениях алгоритма~\ref{alg:ellipsoid})
\[
    \text{diam}_k = 2 \cdot \lambda_{max}^{1/2}(H_k),
\]
где $\lambda_{max}$ обозначает наибольшее собственное значение матрицы.

Другая неточность оракула в методе решения внешней задачи порождается приближенностью решения внутренней задачи, в данном случае быстрым градиентным методом. Пусть внутренний метод настроен на точность $\delta_{FGM}$ и выполняет число итераций, достаточное для достижения этой точности в соответствии с теоретическими оценками. 

Пусть метод решения внешней задачи, то есть метод эллипсоида, работает до выполнения условия останова \eqref{constr_cond}. Будем полагать целью метода достижение им в итоге $\varepsilon$-решения общей $\max$-$\min$ задачи. Тогда, следуя рассуждениям из прошлого примера, необходимо настроить условие останова метода эллипсоидов на точность $\varepsilon - \delta_{FGM} - \delta$. Причём эта величина зависит от способа учёта аддитивной неточности градиента, так же как и практическое время выполнения методом условия останова. Исследуем эту практическую зависимость. 

\begin{figure}[H]
	\centering
	\includegraphics[width=0.7\linewidth]{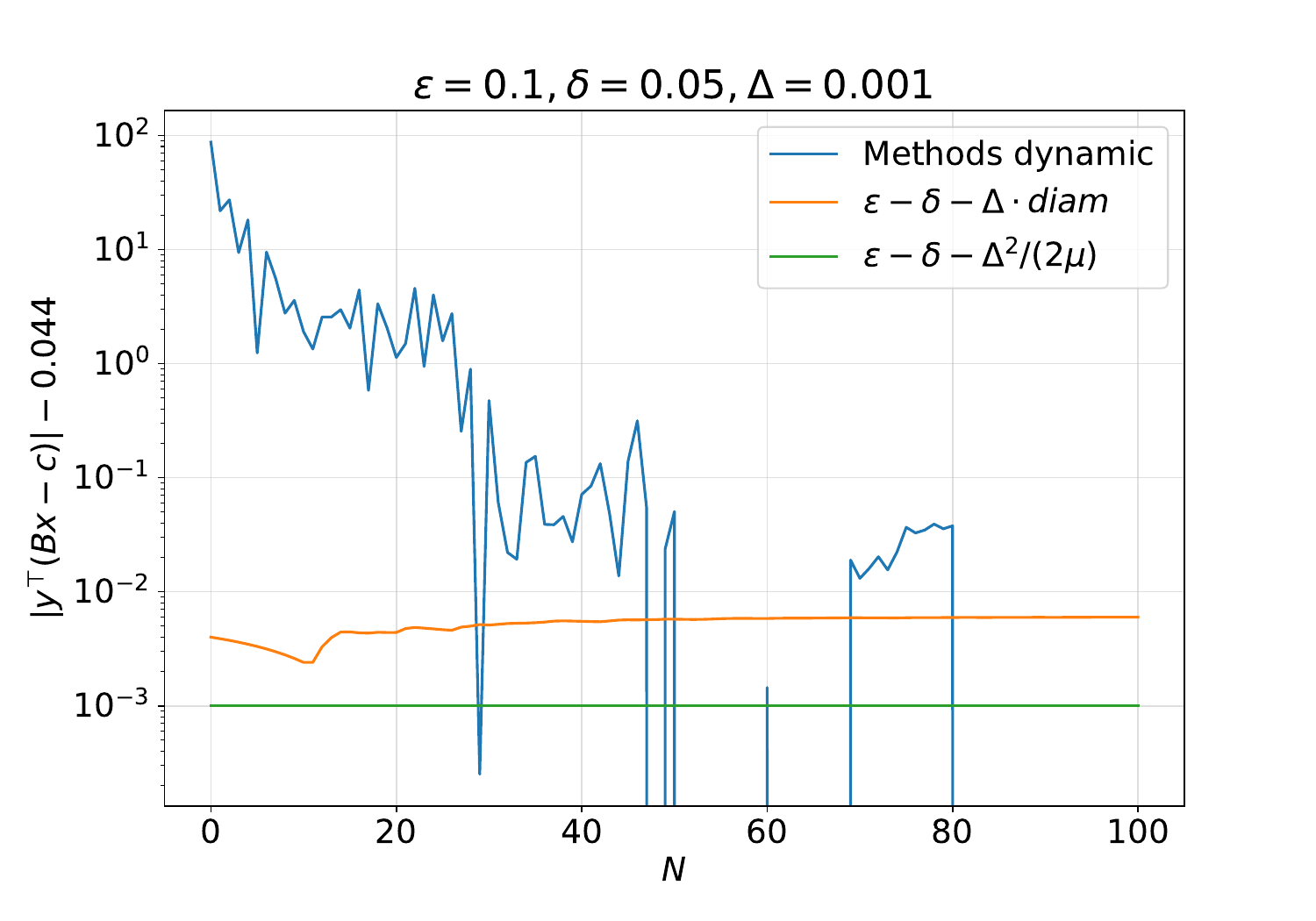}
	\caption{Практическое изменение условия останова метода эллипсоидов для лагранжевой седловой задачи к задаче LogSumExp}
	\label{stopping_condition_evo}
\end{figure}

На рис.~\ref{stopping_condition_evo} представлены графики значений левой части выражения из условия \eqref{constr_cond} и значений $\varepsilon - \delta_{FGM} - \delta$ для двух способов учёта неточности градиента для частного случая задачи при $n=3$, $m=10$, $\varepsilon=0.1$, $\delta=0.05$, $\Delta=0.001$, $\mu = 0.0001$. Как можно видеть, на практике происходит уменьшение диаметра эллипсоидов, и при динамическом учёте неточности градиента граница из условия останова увеличивается заметно в сравнении с динамикой значений, достигаемых в генерируемых методом точках, что может позволить уменьшить число итераций и время работы метода до выполнения этого условия при сохранении гарантий на точность получаемого решения.

\subsection{Задача проектирования точки на множество, определённое набором гладких ограничений}\label{sec:exp-proj}

Теперь рассмотрим задачу проектирования точки на выпуклое множество с нетривиальной структурой, задаваемое набором некоторого небольшого числа ограничений вида неравенств с гладкими функциями \cite{usmanova2021fast}. Проектирование возникает в качестве подзадачи во многих алгоритмах оптимизации, в свою очередь, применяющихся для решения задач с ограничениями. Не всегда проектирование может быть выполнено точно с разумной алгоритмической сложностью. В этой связи разумно ставить задачу нахождения проекции с некоторой точностью, то есть задачу отыскания приближённого решения задачи вида
\begin{align*}
    &\min_{x \in \mathbb{R}^n} \|x_0 - x\|_2^2, \\
    &\text{удовл.} \;\;g_i(x) \leqslant 0,\;\;g_i\text{~---  }L\text{-гладкая},\;\;\forall i=1,...,m,
\end{align*}
лагранжевой седловой задачей к которой является следующая:
\[
    \max_{\lambda \in \mathbb{R}_{+}^m} \min_{x \in \mathbb{R}^n} \left\{ \|x_0 - x\|_2^2 + \sum_{i=1}^m \lambda_i g_i(x)\right\}.
\]
Для случая малого числа ограничений $m$ в работе \cite{usmanova2021fast} был предложен эффективный подход, алгоритмическая сложность которого линейно зависит от размерности $n$, основанный на совместном применении метода эллипсоидов (или метода Вайды) и быстрого градиентного метода. Аналогичный подход, описанный в данной работе, имеет значимые отличия, с одной стороны, в требованиях к точности решения вспомогательной $\min$-задачи (которая в соответствии с предлагаемым анализом может быть выбрана равной $\varepsilon/2$, тогда как в подходе \cite{usmanova2021fast} она необходимо $\sim \varepsilon^4$) и, с другой стороны, в применяемом условии останова \eqref{constr_cond} (гарантирующее достижение заданной точности исходной прямой задачи, которое может оказаться выполненным прежде теоретически достаточного числа итераций, что даёт существенное удобство для приложений). В следующей таблице~\ref{res_proj} указано время выполнения предлагаемого в данной статье подхода (алгоритм~\ref{alg:ellipsoid}, эллипсоиды) и подхода работы \cite{usmanova2021fast} (Algorithm 4, FPM) для случая $m=3$ ограничений вида 
\[g_i(x) = (x - x_i)^\top A_i (x - x_i) - r_i\leqslant 0,\]
где матрицы $A_i$ положительно определённы и генерируются случайно с элементами из $\mathcal{U}(0, 0.05)$, центральные точки $x_i$ имеют независимо случайные компоненты из $\mathcal{U}(-1, 1)$, значения $r_i$ равномерно и независимо случайно сгенерированы из $\mathcal{U}(0, 0.1)$. Проверка итоговой точности осуществляется в сравнении с решением, полученным одним из методов, настроенным на точность $\varepsilon = 10^{-10}$ (как по прямой функции $\|x_0 - x\|^2_2 \leq \|x_0 - x^*\|^2_2 + \varepsilon$, так и по ограничениям $g_i(x) \leq \varepsilon\;\forall i$). Как можно видеть, описанный в этой статье подход за счёт используемого условия останова и уменьшенной трудоёмкости вспомогательных задач на практике оказывается существенно более эффективным в смысле времени выполнения.

\begin{center}
\begin{table}[ht!]
\centering
\begin{tabular}{|c|c|c|c|c|}
\hline
\multirow{3}{*}{$\varepsilon$} & \multicolumn{4}{c|}{Время работы, c} \\ \cline{2-5}
& \multicolumn{2}{c|}{$n=200$} & \multicolumn{2}{c|}{$n=300$} \\ \cline{2-5}
& Эллипсоиды & FPM & Эллипсоиды & FPM \\ \hline

$10^{-1}$ & 2  & 13  & 2  & 15  \\
$10^{-2}$ & 3  & 54  & 12 & 70  \\
$10^{-4}$ & 16 & 119 & 29 & 178 \\
$10^{-5}$ & 30 & 171 & 45 & 271 \\
$10^{-6}$ & 33 & 210 & 47 & 336 \\\hline

\end{tabular}
\caption{Сравнение работы методов для набора из $m = 3$ квадратичных ограничений}
\label{res_proj}
\end{table}
\end{center}


\section{Заключение}\label{section:conclusions}

В настоящей работе получены оценки сложности для сильно выпукло-вогнутых седловых задач вида
\begin{equation}\label{eqconclusion}
    \min_{x \in Q_{x}}\max_{y \in Q_y} \left\{ S(x, y) := r(x) + F(x, y) - h(y) \right\}
\end{equation}
в случае, когда одна из групп переменных ($x$ или $y$) имеет большую размерность, а другая~--- малую (несколько десятков). 

Первые два предлагаемых подхода к такого типа задачам основаны на использовании для подзадачи выпуклой минимизации (вогнутой максимизации) для группы переменных малой размерности методов секущей гиперплоскости (метод эллипсоидов или метод Вайды). Мы приводим оба варианта как с методом эллипсоидов, так и с методом Вайды, поскольку каждый из них имеет свои преимущества: метод Вайды приводит к лучшей оценке числа итераций, а метод эллипсоидов~--- к меньшей сложности итераций по сравнению с методом Вайды. При этом в случае внешней подзадачи небольшой размерности важно использовать эти методы уже в авторском варианте с заменой обычного субградиента на $\delta$-субградиент (отметим, что тут можно использовать и $\delta$-неточный субградиент, причём оценки сложности асимптотически при $\varepsilon \rightarrow 0$ будут теми же). Вспомогательные подзадачи оптимизации по группе переменных большой размерности при этом предлагалось решать с помощью ускоренных градиентных методов. Такая схема позволила вывести приемлемые оценки сложности, зависящие как от обусловленности целевой функции, так и от размерности пространства (см. теорему~\ref{complexity_theorem} и пункт~\ref{second_x}).

Заметим, что первый подход (малая размерность $x$) можно применить и в случае малой размерности $y$, записав аналог задачи \eqref{eqconclusion} следующим образом:
\begin{equation}\label{problem:min_max1}
    \min_{y \in Q_{y}}\max_{x \in Q_{x}} \left\{h(y) - F(x,y) - r(x) \right\}.
\end{equation}

Напомним, что $r$ предполагается проксимально-дружественной функцией, что означает возможность в явном виде решить подзадачу
\begin{equation}\label{prox_r_conclusion}
    \min_{x \in Q_{x}}\left\{\langle c_{1}, x\rangle + r(x)+ c_{2}\|x\|_{2}^{2}\right\},\quad c_1 \in Q_x,\ c_2 > 0.
\end{equation}
В таблице~\ref{comparison} указано количество операций, необходимых для того, чтобы решить задачу \eqref{problem:min_max1} с точностью $\varepsilon$ по $y$ (подход~\ref{appr:first}) или решить аналогичную задачу \eqref{eqconclusion} с точностью $\varepsilon$ по $x$ (подход~\ref{appr:second}).
\begin{table}[htp]
\caption{Сравнение первого и второго подходов}\label{comparison}
\begin{center}
\scalebox{0.88}{
\begin{tabular}{ |c|c|c| }
 \hline
 Подход~\ref{appr:first} & Подход~\ref{appr:second} & Операция \\
 \hline
 $O \left( m \ln \frac{m}{\varepsilon} \right)$ &
 $O \left( m \sqrt{\frac{L_{xx}}{\mu_x} + \frac{2L_{xy}^2}{\mu_x\mu_y}} \ln \frac{m}{\varepsilon} \ln \frac{1}{\varepsilon} \right)$ & вычислений $\nabla_y F,\, \nabla h$ \\
 $O \left( m \sqrt{\frac{L_{xx}}{\mu_x}}  \ln \frac{m}{\varepsilon} \ln \frac{1}{\varepsilon} \right)$ &
 $O \left( \sqrt{\frac{L_{xx}}{\mu_x} + \frac{2L_{xy}^2}{\mu_x\mu_y}} \ln \frac{1}{\varepsilon} \right)$ & вычислений $\nabla_x F$, \eqref{prox_r_conclusion} \\
 \hline
\end{tabular}}
\end{center}
\end{table}
Согласно таблице~\ref{comparison} в большинстве случаев второй подход проигрывает первому. Тем не менее, если $m \ln m \gg \sqrt{\frac{L_{xx}}{\mu_x} + \frac{2L_{xy}^2}{\mu_x\mu_y}}$ и вычисление $\nabla_x F$, \eqref{prox_r_conclusion} является трудоёмким, то второй подход может оказаться эффективнее.

Помимо методов секущей гиперплоскости с неточным $\delta$-субградиентом в работе предложен также некоторый аналог метода дихотомии для решения маломерных задач выпуклой оптимизации с использованием неточных градиентов на итерациях. Этот метод назван {\it многомерной дихотомией}. По сути, он есть обобщение обычной (одномерной) дихотомии на задачи минимизации функций $n$ переменных. Оказалось, что для задач очень малой размерности использование указанного подхода вполне оправдано. Были представлены условия для решения вспомогательной задачи и доказана сходимость метода при выполнении этих условий на каждом шаге. Кроме этого, была получена оценка на достаточное количество итераций для достижения искомой точности по функции (теорема~\ref{th:dich_x}). Данная оценка зависит от размерности пространства (эта зависимость сопоставима с $O\left(2^{n^2}\right)$), а также от требуемой точности решения (эта зависимость имеет вид $O\left(\log_2^n \frac{1}{\varepsilon}\right)$). Полученный результат смотрится значительно хуже по сравнению с описанными оценками для метода эллипсоидов с $\delta$-субградиентом. Однако проведённые численные эксперименты показали, что предложенный метод многомерной дихотомии может работать эффективнее  метода эллипсоидов при $n=2$, что соответствует случаю двух ограничений в прямой задаче.

В ходе экспериментов была сопоставлена работа метода дихотомии, быстрого градиентного метода с $(\delta, L)$-оракулом, а также метода эллипсоидов или метода Вайды с использованием $\delta$-субградиентов для седловых задаче с малой размерностью по одной из переменных. Точнее говоря, выполнены эксперименты на двойственной задаче для задачи минимизации LogSumExp функции с $\ell_2$-регуляризацией в размерности $m$ и с $n$ линейными ограничениями. В ходе этих экспериментов было установлено следующее. Во-первых, маломерные методы быстрее быстрого градиентного метода при высокой требуемой точности. В наших условиях таковой точностью является $\varepsilon=10^{-9}$. Во-вторых, метод многомерной дихотомии быстрее метода эллипсоидов при $n=2$, однако при повышении этой размерности время его работы критически увеличивается, и уже при $n=3$ его эффективность исчезает. В-третьих, было получено, что время работы быстрого градиентного метода для данной задачи увеличивается не так существенно при увеличении $m$, как для метода эллипсоидов или многомерной дихотомии. Кроме того, была сопоставлена работа метода многомерной дихотомии, рассматриваемых в работе методов секущей гиперплоскости (методы эллипсоидов и Вайды), а также быстрого градиентного метода для задачи минимизации квадратичной функции (для $n=400$ и $n=1000$) с двумя негладкими ограничениями, $\max$-агрегирующими несколько ($m=10, 20$) линейных ограничений. Для данной задачи время работы метода дихотомии и его вариантов (метод на треугольнике) оказалось меньше времени работы метода эллипсоидов (как для всех исходных, так и для агрегированных ограничений) и быстрого градиентного метода. Было произведено сравнение способов учёта неточностей, возникающих при наличии аддитивного шума в значениях градиента, для случая применения к задаче малой размерности метода эллипсоидов. При значении неточности, изменяющемся вместе с диаметром текущего эллипсоида, применяемый метод достигает выполнения условия останова быстрее, чем при постоянной оценке неточности. Сравнение работы методов было произведено также для задачи проектирования точки на множество, заданное набором гладких функционалов ограничений \cite{usmanova2021fast}. Предлагаемый нами  подход применения метода эллипсоидов для подзадач небольшой размерности при использовании предложенного условия останова оказывается эффективнее по сравнению с алгоритмом, предложенным в работе \cite{usmanova2021fast} для аналогичной постановки задачи. В рамках поставленных численных экспериментов рассматривались постановки задач, для которых задача малой размерности (в данном случае по двойственным переменным лагранжевой седловой задачи) не является сильно выпуклой (вогнутой). Несмотря на то что теоретический анализ оценок скорости методов в настоящей статье приводится для случая сильно выпукло-вогнутых задач, правильная настройка методов на практике позволяет с тем же успехом применять предложенные схемы и для просто выпуклых (или вогнутых маломерных подзадач с обеспечением достижения желаемого качества решения задачи. Это объясняется отсутствием необходимости требовать сильную выпуклость (вогнутость) целевой функции (она важна лишь для теоретических оценок) для реализации всех применяемых в работе к подзадачам небольшой размерности методов.

\end{fulltext}

\appendix

\section{Доказательство леммы~\ref{lem:ellips}}
Далее, приводится доказательство из \cite{polyak1983intro} (с. 123--124),
где вместо сильной вогнутости $S$ по $y$ используется предположение о компактности $Q_y$.

Пусть $\nu \in \partial_x S(x, \widetilde{y})$. Для любого $x' \in Q_x$
$$
    \widehat{g}(x') = \max_{y \in Q_y} S(x', y) \geq S(x', \widetilde{y}) \geq S(x, \widetilde{y}) + \langle \nu, x' - x \rangle \geq \widehat{g}(x) + \langle \nu, x' - x \rangle - \delta.
$$
Таким образом, $\nu \in \partial_{\delta} \widehat{g}(x)$, что и требовалось доказать.

\section{Доказательство теоремы~\ref{lem:boundary}}
Заметим, что $(\delta,L)$-субградиент $\nabla_{\delta, L} g(x)$ в определении~\ref{deltaLsubgrad} соответствует $(2\delta,L)$-оракулу вида $(g(x)-\delta, \nabla_{\delta, L} g(x))$ в определении 1 из работы \cite{th2_cite}, то есть выполнено следующее неравенство:
\begin{equation}
0 \leqslant g(y) -\left(g(x)-\delta + \langle \nabla_{\delta, L} g(x), y - x \rangle\right)\leqslant \frac{L}{2}\|y - x\|_2^2 + 2\delta.
\end{equation}
В разделе 2.2 в \cite{th2_cite} было показано, что если $\rho(x,\partial Q_x)\geq 2\sqrt{\frac{\delta}{L}}$, то для любого субградиента $\nabla g(x)$ верно следующее утверждение: $$\|\nabla_{\delta, L} g(x)-\nabla g(x)\|\leq 2\sqrt{\delta L},$$
что и требовалось доказать.


\section{Доказательство леммы~\ref{lem:delta_vs_L}}
Из сильной вогнутости $S$ по $y$ следует, что для любого $x$ задача максимизации \eqref{problem:max_S0} имеет единственное решение, которое мы будем обозначать за $y(x)$, а также то, что для любого $y \in Q_y$ выполняется неравенство
\begin{equation*}
    S(x, y) \leq \underbrace{S(x, y(x))}_{\widehat{g}(x)} - \frac{\mu_y}{2} \|y - y(x) \|_2^2.
\end{equation*}
В частности, если $\widetilde{y}$~--- $\tilde{\varepsilon}$-решение внутренней задачи \eqref{problem:max_S0}, то
\begin{equation}\label{arg_convergence}
    \|\widetilde{y} - y(x) \|_2^2 \leq \frac{2}{\mu_y} \tilde{\varepsilon}.
\end{equation}
Согласно теореме Демьянова-Данскина \cite{bernhard1995theorem,DDR-theorem} в любой точке $x \in Q_x$ функция $\widehat{g}$ дифференцируема, и её градиент равен
\begin{equation}\label{dansk}
    \nabla \widehat{g}(x) = \nabla_x S(x, y(x)).
\end{equation}
Используя \eqref{Lxy0}, \eqref{arg_convergence} и \eqref{dansk}, получим
\begin{equation*}
    \|\nabla_x S (x, \widetilde{y}) - \nabla \widehat{g}(x)\|_2 \leqslant L_{xy} \sqrt{\frac{2 \tilde{\varepsilon}}{\mu_y}},
\end{equation*}
что и требовалось доказать.

\section{Доказательство леммы~\ref{delta_12}}
\begin{enumerate}
    \item 
Для любых $x,x' \in Q,\, \nabla g(x) \in \partial g(x)$ и $\delta_1$-неточного субградиента $\nu$ в точке $x$
\begin{equation*}
    \begin{aligned}
    g(x') &\geqslant g(x) + \left\langle \nabla g(x), x'-x  \right\rangle \\&
     = g(x) + \left\langle \nu, x'-x \right\rangle + \left\langle \nabla g(x) -\nu, x'-x \right\rangle \\&
     \geqslant g(x) + \left\langle \nu, x'-x \right\rangle - \delta_1 \operatorname{diam} Q.
\end{aligned}
\end{equation*}
Таким образом, $\delta_1$-неточный градиент $\nu$ является $\delta_2$-субградиентом $g$ в точке $x$ с $\delta_2 = \delta_1 \operatorname{diam} Q$.
    \item 
Для любых $x,x' \in Q,\, \nabla g(x) \in \partial g(x)$ и $\delta_1$-неточного субградиента $\nu$ в точке $x$
\begin{equation*}
    \begin{aligned}
    g(x') &\geqslant g(x) + \left\langle \nabla g(x), x'-x  \right\rangle +\frac{\mu}{2}\|x'-x\|_2^2 \\&
     = g(x) + \left\langle \nu, x'-x \right\rangle + \left\langle \nabla g(x) -\nu, x'-x \right\rangle +\frac{\mu}{2}\|x'-x\|_2^2 \\&
     \geqslant g(x) + \left\langle \nu, x'-x \right\rangle - \delta_1 \|x'-x\|_2 +\frac{\mu}{2}\|x'-x\|_2^2.
\end{aligned}
\end{equation*}
Учитывая, что $\delta_1 \|x'-x\|_2 \leqslant \frac{\mu}{2}\|x'-x\|_2^2 + \frac{\delta_1^2}{2\mu}$, получим
$$
g(x') \geqslant g(x) + \left\langle \nu, x'-x \right\rangle -  \frac{\delta_1^2}{2\mu}.
$$
Таким образом, $\delta_1$-неточный градиент $\nu$ является $\delta_2$-субградиентом $g$ в точке $x$ с $\delta_2 = \frac{\delta_1^2}{2\mu}$.
\end{enumerate}

\section{Доказательство леммы~\ref{subgradient}}


Если $\textbf{x}_*$ есть внутренняя точка, то градиент по нефиксированным переменным равен нулю в силу того, что $\textbf{x}_*$~--- минимум. Тогда с учетом того, что $\nabla f(\textbf{x}_*)\in\partial f(\textbf{x}_*)$, получаем утверждение леммы.

Допустим, что $\textbf{x}_*$~--- граничная точка. Тогда множество условного субдифференциала на гиперкубе $Q$ определяется следующим образом:
$$\partial_Q f(\textbf{x}) = \partial f(\textbf{x}) + N\left(\textbf{x}\;|\;Q\right),$$
где $N(\textbf{x}\;|\;Q) = \left\{\textbf{a}\;|\;\langle\textbf{a}, \textbf{y} - \textbf{x}\rangle\leqslant 0\; \forall \textbf{y} \in Q\right\}$.

В случае дифференцируемой функции имеем:
$$\partial f(\textbf{x}_*) = \{\nabla f(\textbf{x}_*)\}.$$

Из того, что $\textbf{x}_*$ есть 
граничная
точка, следует, что существует непустой набор координат $\{x_j\}_j$ такой, что $x_j = \max_{\textbf{y}\in Q_k}y_j$ или $x_j = \min_{\textbf{y}\in Q_k}y_j$. Введём обозначения:
$$
    J_+=\{j\in \mathbb{N}\;|\; x_j = \max_{\textbf{y}\in Q_k}y_j\},
$$
$$
    J_-=\{j\in \mathbb{N}\;|\;x_j = \min_{\textbf{y}\in Q_k}y_j\}.
$$

В таком случае заметим, что любой вектор $a$ такой, что $a_j \geqslant 0\;\forall j \in J_+$, $a_j \leqslant 0\;\forall j \in J_-$ и $a_j=0$ в противном случае, принадлежит нормальному конусу.

Также заметим, что $(\nabla f(\textbf{x}_*))_j \leqslant 0\;\forall j \in J_+$, $(\nabla f(\textbf{x}_*))_j\geqslant 0\;\forall j \in J_-$ и $(\nabla f(\textbf{x}_*))_j=0$ в противном случае. Действительно, если $(\nabla f(\textbf{x}_*))_j > 0$ для некоторого $j\in J_+$, то существует вектор $\textbf{x}=\textbf{x}_* + \alpha \textbf{e}_k \in Q$ для некоторого $\alpha<0$ и вектора $e_k^j = \delta_{kj}$, где $\delta_{kj}=1$ для $k=j$ и $\delta_{kj}=0$ иначе. Причём значение функции в этой точке будет $f(\textbf{x}) = f(\textbf{x}_*) + \alpha (\nabla f(\textbf{x}_*))_j + o(\alpha) < f(\textbf{x}_*)$ для достаточно малого $\alpha$, что противоречит тому, что $\textbf{x}_*$~--- решение.

Тогда, выбрав $\textbf{a}$ такой, что $\textbf{a}_\parallel = -\left(\nabla f(\textbf{x}_*)\right)_\parallel$, получаем субградиент из условия:
$$
    \textbf{g} = \nabla f(\textbf{x}_*) + \textbf{a},\;\; \textbf{g}_\parallel =0.
$$

\section{Доказательство теоремы~\ref{th_ellips}}

В нашем методе, как и в обычном методе эллипсоидов, на каждом шаге эллипсоид рассекается плоскостью, проходящей через его центр, а затем рассматривается эллипсоид наименьшего объёма, содержащий одну из частей. Можно доказать  (см., например, \cite{bubeck2015convex}), что на каждом шаге выполняется неравенство:
\begin{equation}\label{volume}
    \frac{vol (\mathcal{E}_{k+1})}{vol(\mathcal{E}_k)} \leqslant e^{-1/2n}\ \Rightarrow\ vol(\mathcal{E}_N) \leqslant e^{-N/2n} vol(\mathcal{B}_{\mathcal{R}}). 
\end{equation}
Если $w_k = 0$, то по определению $\delta$-субградиента $g(x) \geqslant g(c_k) - \delta\ \forall x \in Q_x\ \Rightarrow g(c_k) - g(x_*) \leqslant \delta$, и условие теоремы выполнено. Далее, считаем, что вектор $w_k$ ненулевой. Если $c_k \in Q_x$, то в силу определения $\delta$-субградиента справедливо включение
\begin{equation} \label{remove}
    (\mathcal{E}_k \setminus \mathcal{E}_{k+1}) \cap Q_x \subseteq \{ x \in Q_x : \langle w_k, x - c_k \rangle > 0 \} \subseteq \{ x \in Q_x : g(x) > g(c_k) - \delta \}.
\end{equation}
Для $\varepsilon \in [0,1] $ рассмотрим множество $Q_x^{\varepsilon} := \{(1 - \varepsilon)x_* + \varepsilon x, x \in Q_x \}$. Заметим, что $Q_x^{\varepsilon} \subseteq \mathcal{E}_0$ и
$$
    vol(Q_x^{\varepsilon}) = \varepsilon^n vol(Q_x) \geqslant \varepsilon^n vol(\mathcal{B}_{\rho}) = \left(\frac{\varepsilon \rho}{\mathcal{R}} \right)^n vol(\mathcal{B}_{\mathcal{R}}).
$$
Для $\varepsilon > e^{-N/2n^2}\frac{\mathcal{R}}{\rho}$ получаем, что из \eqref{volume} следует, что $vol(Q_x^{\varepsilon}) > vol(\mathcal{E}_N) \Rightarrow $ найдётся шаг $j \in \{0, \ldots, N-1 \}$, и найдутся $x_{\varepsilon} \in Q_x^{\varepsilon}$ такие, что $x_{\varepsilon} \in \mathcal{E}_j$ и $x_{\varepsilon} \notin \mathcal{E}_{j+1}$. Если бы при этом точка $c_j$ лежала вне $Q_x$, то мы бы отсекли часть эллипсоида $\mathcal{E}_j$, не пересекающуюся с $Q_x$, что привело бы к противоречию с тем, что $x_{\varepsilon} \in Q_x$. Значит, $c_j \in Q_x$. Тогда, воспользовавшись \eqref{remove}, получим, что $g(x_{\varepsilon}) > g(c_j) - \delta$. Поскольку $\exists x \in Q_x: x_{\varepsilon} = (1-\varepsilon)x_*+\varepsilon x$, то в силу выпуклости $g$ имеем
$$
    g(x_{\varepsilon}) \leqslant (1-\varepsilon) g(x_*)+\varepsilon g(x) \leqslant (1-\varepsilon) g(x_*)+\varepsilon \left((g(x_*) + B \right) = g(x_*) + B\varepsilon.
$$

Получаем, что
\begin{equation}\label{ellip_proof}
\begin{split}
    g(c_j) < g(x_*) + B\varepsilon + \delta\ \forall \varepsilon > e^{-N/2n^2}\frac{\mathcal{R}}{\rho}\ \Rightarrow \\
    \Rightarrow g(c_j) - g(x_*) \leqslant e^{-N/2n^2}\frac{B\mathcal{R}}{\rho} + \delta,
\end{split}
\end{equation}
откуда следует \eqref{th_ellips_1}. Если дополнительно $g$ является $\mu$-сильно выпуклой, то есть
$$
    g(x) - g(x') - \langle \nabla g(x'), x - x' \rangle \geqslant \frac{\mu}{2} \| x - x' \|_2^2\quad \forall x, x' \in Q_x,
$$
то, подставив $x = c_j,\ x' = x_*$ и использовав $\langle \nabla g(x_*), x - x_*' \rangle\ \forall x \in Q_x$, получим
$$
    g(c_j) - g(x_*) \geqslant \frac{\mu}{2} \| c_j - x_* \|_2^2.
$$
Отсюда, а также и из \eqref{ellip_proof} следует второе доказываемое утверждение \eqref{th_ellips_2}.

\section{Доказательство теоремы~\ref{th:dich_x}}

Пусть мы решаем задачу вида
\begin{equation}
\label{func_f2}
\min_x f(x).
\end{equation}

Оценим сложность метода многомерной дихотомии (то есть количество обращений к подпрограмме вычисления градиента $\nabla f$, достаточное для достижения $\varepsilon$-точного решения по функции).

В доказательстве данной теоремы будем использовать обоснованную в теореме~\ref{FullCond} оценку необходимого количества внешних итераций для достижения приемлемого качества приближённого решения задачи минимизации $f$:
\begin{equation}
N = \left\lceil\log_2 \left(\frac{4R(M_f+2L_f R)}{L_f\varepsilon}\right)\right\rceil.
\end{equation}


Пусть $T(n, R, \varepsilon)$~--- количество вспомогательных задач минимизации соответствующей функции размерности $n-1$, которых достаточно для решения задачи в размерности $n$ на гиперкубе диаметра $R$ и точностью $\varepsilon$. Для $n=0$ положим $T(0, R, \varepsilon)$. Заметим, что одна итерация требует решения $n$ вспомогательных задач. С учетом этого получаем следующую
рекуррентную формулу для основной задачи:
$$
T(n,R,\varepsilon)=\sum\limits_{k=0}^{\left\lceil\log_2 \left(\frac{M_f R}{\varepsilon}\right)\right\rceil}n T(n-1,R \cdot 2^{-k},\widetilde{\varepsilon})
$$
и аналогичное выражение уже с учётом всех необходимых вспомогательных подзадач:
$$T(n,R,\varepsilon)=\sum\limits_{k=0}^{\left\lceil\log_2 \left(\frac{C_1 R}{\varepsilon}\right)\right\rceil}n T(n-1,R\cdot 2^{-k}, \widetilde{\varepsilon}).$$
где $\widetilde{\varepsilon}$ определяется согласно \eqref{estimate_gen_dih}, а $$C_1 = \max\left(M_f, \frac{4(M_f +2L_f R)}{L_f}\right),$$ 
Пусть $C_\varepsilon = \frac{128 L_f^2}{\mu_f}$. Докажем индукцией по $n$ следующую оценку:
\begin{equation}
\label{RecIneq}
T(n, R, \varepsilon) \leqslant 2^{\frac{n^2+n}{2}}\log^n_2 \left(\frac{CR}{\varepsilon}\right)+O\left(\log_2^n \left(\frac{C R}{\varepsilon}\right)\right), \text{ где } C = 2\max{\left(C_1, C_\varepsilon\right)}.
\end{equation}
Во введенных выше обозначениях коэффициент 2 в выражении для $C$ позволяет избежать записи округления вверх в дальнейшем.

Базис индукции очевиден: 
$$
    T(1, R, \varepsilon) = \log_2\frac{C_1}{\varepsilon} \leqslant \log_2\left(\frac{CR}{\varepsilon}\right).
$$

Допустим справедливость \eqref{RecIneq} для некоторой размерности $n$ и докажем, что \eqref{RecIneq} верно и для размерности $n+1$.
$$
    T(n+1,R,\varepsilon)=\sum\limits_{k=0}^{\left\lceil \log_2\left(\frac{C_1 R}{\varepsilon}\right)\right\rceil}(n+1)T(n,R\cdot2^{-k}, \widetilde{\varepsilon})\leq
$$
$$
    \leq (n+1) \cdot 2^{\frac{n^2+n}{2}} \sum\limits_{k=0}^{\left\lceil\log_2 \left(\frac{C_1 R}{\varepsilon}\right)\right\rceil} \log^n_2 \left(\frac{CC_\varepsilon R^2}{2^{2k} \varepsilon^2}\right)+O\left(\log_2^n \left(\frac{C R}{\varepsilon}\right)\right).
$$
Оценим сумму:
\begin{equation*}
    \begin{aligned}
    \sum\limits_{k=0}^{\left\lceil\log_2 \left(\frac{C_1 R}{\varepsilon}\right)\right\rceil} \log^n_2 \left(\frac{CC_\varepsilon R^2}{2^{2k} \varepsilon^2}\right) & \leq \sum\limits_{k=0}^{\left\lceil\log_2 \left(\frac{C R}{\varepsilon}\right)\right\rceil} \log^n_2 \left(\frac{C^2 R^2}{2^{2k} \varepsilon^2}\right) 
    \\ & \leq 2^n \cdot \int\limits_{0}^{\log_2 \left(\frac{C R}{\varepsilon}\right)+1}\left(\log_2 \left(\frac{C R}{\varepsilon}\right) - k\right)^n dk+\log_2^n \left(\frac{C R}{\varepsilon}\right)
    \\& = \frac{2^n}{n+1}\left(\log_2^{n+1} \left(\frac{C R}{\varepsilon}\right)+1\right)+\log_2^n \left(\frac{C R}{\varepsilon}\right).
    \end{aligned}
\end{equation*}
Поэтому верно неравенство
$$
    T(n+1,R,\varepsilon) \leq 2^{\frac{(n+1)^2+(n+1)}{2}} \log_2^{n+1} \left(\frac{C R}{\varepsilon}\right)+O\left(\log_2^n \left(\frac{C R}{\varepsilon}\right)\right),
$$
откуда получаем доказываемую оценку \eqref{RecIneq}. Окончательно получаем, что для решения задачи \eqref{func_f2} достаточно следующего числа
\begin{equation}
O\left(2^{n^2} \log_2^n \left(\frac{C R}{\varepsilon}\right)\right),\text{ где } C = \max\left(M_f,\frac{4(M_f+2L_fR)}{L_f}, \frac{128L_f^2}{\mu_f}\right)
\end{equation}
вычислений неточного градиента $\nu(\textbf{x})$. При этом любую вспомогательную задачу для текущего уровня рекурсии мы решаем с точностью по аргументу (см. \eqref{tilde_delta})
$$\widetilde{\delta}=\frac{\Delta}{C_f}\geqslant 2^{-N},$$
где $N$ определяется в утверждении теоремы~\ref{FullCond}.

\section{ Доказательство теоремы~\ref{InexGradConst}}

В работе \cite{Ston_Pas} было доказано, что если задан гиперкуб $Q$ с максимальным расстоянием между точками $R$ и требуется минимизировать на нём функцию с точностью $\varepsilon$, то для этого в методе достаточно решить вспомогательную задачу с точностью по аргументу
\begin{equation}
    \Delta \leqslant \frac{\varepsilon}{8L_f R},
\end{equation}
где $R$ есть размер начального гиперкуба. Данная теорема доказана для размерности $n=2$ в работе \cite{Ston_Pas}, однако она без труда обобщается на большие размерности. Если $f$ есть $\mu_f$-сильно выпуклая функция, то, используя написанное выше условие останова, получаем, что достаточно решить вспомогательную задачу с точностью $\widetilde{\varepsilon}$ по функции:
\begin{equation}
    \widetilde{\varepsilon} \leqslant \frac{\mu_f \varepsilon^2}{128 L_f^2 R^2}.
\end{equation}

\section{ Доказательство теоремы~\ref{CurGrad}}



Далее, нам понадобится следующее очевидное соотношение:
\begin{equation}\label{trick_lemma}
\forall a,b\in\mathbb{R}\;\; |a-b|\leqslant |b| \Rightarrow ab \geqslant 0.
\end{equation}

Заметим, что множество $Q_k$ на $k$-ой итерации выбирается правильно, если знак производной в решении вспомогательной задачи по фиксированной переменной совпадает со знаком производной в приближении решения.

Пусть $\nu(\textbf{x}) = \nabla f(\textbf{x})$. Из \eqref{trick_lemma} следует, что для того, чтобы совпали знаки $\nu_{\perp Q_k}(\textbf{x}_*)$ и $\nu_{\perp Q_k}(\textbf{x})$, достаточно потребовать
$$\left|\nu_{\perp Q_k}(\textbf{x}_*) - \nu_{\perp Q_k}(\textbf{x})\right| \leqslant |\nu_{\perp Q_k}(\textbf{x})|,$$
где $\nu_{\perp Q_k}(\textbf{x})$~--- проекция вектора $\nu(\textbf{x})$ на ортогональное дополнение к множеству, на котором решается вспомогательная задача.

Используя липшицевость градиента целевого функционала $f$, получаем утверждение теоремы.

\section{ Доказательство теоремы~\ref{small}}

Из леммы~\ref{subgradient} имеем:
$$\textbf{g} \in \partial_Q f(\textbf{x}_*): \textbf{g}_\parallel = 0.$$

Тогда по определению субградиента $f$ в точке $\textbf{x}_*$:
$$f(\textbf{x}^*) - f(\textbf{x}_*) \geq \langle\textbf{g}, \textbf{x}^* - \textbf{x}_*\rangle.$$

Используем неравенство Коши--Буняковского--Шварца:
$$f(\textbf{x}_*) - f(\textbf{x}^*) \leq \|\textbf{g}\|_2a\sqrt{n}.$$

С другой стороны, из условия Липшица $f$ для любой точки $\textbf{x}$ из\\$\Delta$-окрестности точки $
\textbf{x}_*$ имеем:
$$f(\textbf{x})-f(\textbf{x}_*)\leq M_f \Delta,$$
$$f(\textbf{x})-f(\textbf{x}^*) \leq M_f \Delta +\|\textbf{g}\|_2a\sqrt{n} = M_f \Delta + |\nu_{\perp Q_k}(\textbf{x}_*)|R.$$

Ввидy липшицевости градиента $f$ имеем:
$$f(\textbf{x})-f(\textbf{x}^*) \leq M_f\Delta + \left(|\nu_{\perp Q_k}(\textbf{x})|+L_f\Delta\right)R.$$

Пусть после 11-15 шагов Алгоритма \ref{alg:Dichotomy} осталось множество с диаметром $\Delta$ во вспомогательной задаче. Тогда для достижения точности $\varepsilon$ по функции в исходной задаче в некоторой точке $\textbf{x}$ из этого множества достаточно следующего условия:

$$M_f\Delta +\|\textbf{g}\|_2a\sqrt{n} = M_f\Delta + \left(|f_\perp'(\textbf{x})|+L_f\Delta\right)R \leq \varepsilon,$$

$$\Delta\left(M_f + L_f R\right) \leq \varepsilon - |f_\perp'(\textbf{x})|R.$$

Окончательно получаем:
$$\Delta \leq \frac{\varepsilon - R |f_\perp'(\textbf{x})|}{M_f+L_f R}.$$

\section{Доказательство теоремы~\ref{FullCond}}

Объединяя оценки из теорем~\ref{CurGrad} и~\ref{small}, получаем, что для того, чтобы достигнуть точность $\varepsilon$ по функции при решении задачи минимизации на гиперкубе $Q$, каждую вспомогательную задачу нужно решать до тех пор, пока не будет выполнено следующее условие на расстояние от приближенного решения до истинного решения этой задачи:

\begin{equation}
\label{Adaptive}
\Delta \leqslant \max\left\{
	\frac{|\nu_{\perp Q_k}(\textbf{x})|}{L_f},
	\frac{\varepsilon - R |\nu_{\perp Q_k}(\textbf{x})|}{M_f+L_f R}
	\right\}.
\end{equation}

Это условие верно для $\nu(\textbf{x}) = \nabla f(\textbf{x})$. Пусть $\nu(\textbf{x})$ есть такой вектор, что
\begin{equation}
\|\nabla f(\textbf{x}) - \nu(\textbf{x})\|_2 \leqslant \widetilde{\delta}(\textbf{x}).
\end{equation}
В таком случае очевидно, что условие \eqref{Adaptive} будет выполнено, если
\begin{equation}
C_f\widetilde{\delta}(\textbf{x}) + \Delta \leqslant \max\left\{
	\frac{|\nu_{\perp Q_k}(\textbf{x})|}{L_f},
	\frac{\varepsilon - R |\nu_{\perp Q_k}(\textbf{x})|}{M_f+L_f R}
	\right\},
\end{equation}
где $C_f = \max\left(\frac{1}{L_f}, \frac{R}{M_f+L_f R}\right).$

Оценим необходимое число итераций. Если решать вспомогательную задачу текущего уровня рекурсии с точностью $\widetilde{\delta} = \frac{1}{C_f}\Delta,$ то получаем следующее условие:
\begin{equation}
\Delta \leqslant \frac{1}{2}\max\left\{
	\frac{|\nu_{\perp Q_k}(\textbf{x})|}{L_f},
	\frac{\varepsilon - R |\nu_{\perp Q_k}(\textbf{x})|}{M_f+L_f R}
	\right\}.
\end{equation}

Пусть модуль ортогональной компоненты приближения градиента $|\nu_{\perp Q_k}(\textbf{x}_*)|$ равен $q$. Обозначим её приближение с точностью $\Delta$ по аргументу  $\textbf{x}_\Delta$.

После $N$ итераций метода многомерной дихотомии для вспомогательной задачи, используя липшицевость градиента с константой $L_f$, можем получить следующую оценку градиента в точке $\textbf{x}_\Delta$:
$$
    q-2L_f R\cdot 2^{-N} \leqslant |\nu_{\perp Q_k}(\textbf{x})|\leqslant q+2L_f R\cdot 2^{-N}.
$$
В выше написанном неравенстве учитывалось, что размер множества уменьшится в $2^{N}$ раз после $N$ итераций дихотомии, то есть $\Delta\leq 2^{-N}R$ после $N$ итераций.

Значит, для выполнения условия $\Delta \leq \frac{1}{2}
	\frac{|\nu_{\perp Q_k}(\textbf{x})|}{L_f}$ достаточно выполнения следующего условия:
$$R\cdot 2^{-N} \leqslant \frac{q-2L_f R\cdot 2^{-N}}{2L_f }.$$

Аналогичное условие получаем для второй альтернативы $\Delta \leqslant \frac{1}{2}
	\frac{\varepsilon - R |\nu_{\perp Q_k}(\textbf{x})|}{M_f+L_f R}$:
$$
    R\cdot 2^{-N} \leqslant \frac{1}{2}\frac{\varepsilon - qR}{M_f+L_f R} - R\cdot 2^{-N}.
$$

Тогда получаем следующую оценку на $N$:
$$
    R\cdot 2^{-N} \leqslant \frac{1}{4}\min_{q\geqslant 0}\max\left(\frac{q}{L_f}, \frac{\varepsilon - q R}{M_f+L_f R}\right)=\frac{1}{4}\frac{L_f}{M_f+2L_f R}\cdot\varepsilon.
$$

Таким образом, количество итераций исходного алгоритма, необходимое для достижения требуемой точности во вспомогательных задачах, не превосходит следующей величины:
\begin{equation}
\label{dich_N}
N = \left\lceil\log_2 \left(\frac{4R(M_f+2L_f R)}{L_f\varepsilon}\right)\right\rceil.
\end{equation}

\section{ Формулировки некоторых известных используемых результатов}


\begin{theorem}\label{constrains_y}(\cite{gasnikov2018mpt}, упражнение 4.1)
Рассмотрим задачу
$$
    \min_{x\in\mathbb{R}^m} f(x), \quad \text{удовл.}\;\; g(x)\leqslant 0, \, g : \mathbb{R}^m \rightarrow \mathbb{R}^n,
$$
где $f$ и $g_i$ --- выпуклые функции. Лагранжиан этой задачи имеет вид
$$
    \phi(y) = \min_x \left\{f(x)+y^\top g(x)\right\}.
$$
Пусть $x_0$ есть такая точка, что $g(x_0) < 0.$ Тогда для решения $y_*$ задачи $\max_y \phi(y)$ верно следующее неравенство:
\begin{equation}
    \|y_*\|_2 \leqslant \frac{1}{\gamma}\left(f(x_0) - \min_x f(x)\right),
\end{equation}
где $\gamma = \min_k \{-g_k(x_0)\}$.
\end{theorem}

\begin{theorem}\label{dual_L_from_mu}Рассмотрим задачу
$$
    \min_{x \in \mathbb{R}^m} f(x) \quad \text{при}\;\; g(x)\leq 0 \in\mathbb{R}^n.
$$

Пусть $f$ есть $\mu_f$-сильно выпуклая функция, а вектор-функция $g$ есть $M_g$-лип\-ши\-це\-ва функция. Тогда функция $\phi(y) = \min_x \left(f(x)+y^\top g(x)\right)$ имеет липшицев градиент с константой $L = \frac{M_g^2}{\mu_f}$.
\end{theorem}
Доказательство этой теоремы приведено в \cite{Stonyakin} для случая, когда $g(x)$ одномерна, то есть для случая одного условия. Обобщенный результат, представленный здесь, доказывается аналогично.



\end{document}